\documentclass[reqno]{amsart}

\usepackage[all]{xy}

\def\E{\ifmmode{\mathbb E}\else{$\mathbb E$}\fi} 
\def\N{\ifmmode{\mathbb N}\else{$\mathbb N$}\fi} 
\def\R{\ifmmode{\mathbb R}\else{$\mathbb R$}\fi} 
\def\Q{\ifmmode{\mathbb Q}\else{$\mathbb Q$}\fi} 
\def\C{\ifmmode{\mathbb C}\else{$\mathbb C$}\fi} 
\def\H{\ifmmode{\mathbb H}\else{$\mathbb H$}\fi} 
\def\Z{\ifmmode{\mathbb Z}\else{$\mathbb Z$}\fi} 
\def\P{\ifmmode{\mathbb P}\else{$\mathbb P$}\fi} 
\def\T{\ifmmode{\mathbb T}\else{$\mathbb T$}\fi} 
\def\SS{\ifmmode{\mathbb S}\else{$\mathbb S$}\fi} 
\def\DD{\ifmmode{\mathbb D}\else{$\mathbb D$}\fi} 

\newcommand{\e}{\varepsilon}

\newcommand{\del}{\partial}

\newcommand{\ben}{\begin{enumerate}}
\newcommand{\een}{\end{enumerate}}
\newcommand{\be}{\begin{equation}}
\newcommand{\ee}{\end{equation}}
\newcommand{\bea}{\begin{eqnarray}}
\newcommand{\eea}{\end{eqnarray}}
\newcommand{\beastar}{\begin{eqnarray*}}
\newcommand{\eeastar}{\end{eqnarray*}}
\newcommand{\bc}{\begin{center}}
\newcommand{\ec}{\end{center}}

\theoremstyle{theorem}
\newtheorem{thm}{Theorem}[section]

\newtheorem{lem}[thm]{Lemma}
\newtheorem{prop}[thm]{Proposition}

\theoremstyle{definition}
\newtheorem{defn}{Definition}[section]
\newtheorem{rem}[defn]{Remark}

\numberwithin{equation}{section}

\def\R{{\mathbb R}}

\def\E{{\mathbb E}}
\def\Z{{\mathbb Z}}
\def\C{{\mathbb C}}
\def\R{{\mathbb R}}
\def\P{{\mathbb P}}

\def\N{{\mathbb N}}

\def\11{{\mathbb I}}

\def\C{\mathbb{C}}
\def\Z{\mathbb{Z}}

\def\T{\mathbb{T}}

\def\Q{\mathbb{Q}}

\def\E{\ifmmode{\mathbb E}\else{$\mathbb E$}\fi} 
\def\N{\ifmmode{\mathbb N}\else{$\mathbb N$}\fi} 
\def\R{\ifmmode{\mathbb R}\else{$\mathbb R$}\fi} 
\def\Q{\ifmmode{\mathbb Q}\else{$\mathbb Q$}\fi} 
\def\C{\ifmmode{\mathbb C}\else{$\mathbb C$}\fi} 
\def\H{\ifmmode{\mathbb H}\else{$\mathbb H$}\fi} 
\def\Z{\ifmmode{\mathbb Z}\else{$\mathbb Z$}\fi} 
\def\P{\ifmmode{\mathbb P}\else{$\mathbb P$}\fi} 
\def\SS{\ifmmode{\mathbb S}\else{$\mathbb S$}\fi} 
\def\DD{\ifmmode{\mathbb D}\else{$\mathbb D$}\fi} 

\def\R{{\mathbb R}}

\def\E{{\mathbb E}}
\def\Z{{\mathbb Z}}
\def\C{{\mathbb C}}
\def\R{{\mathbb R}}

\def\N{{\mathbb N}}
\def\LL{{\mathcal L}}
\def\MM{{\mathcal M}}

\def\FF{{\mathcal F}}
\def\UU{{\mathcal U}}

\def\e{\varepsilon}

\def\CM{{\mathcal M}}

\def\darr#1{\raise1.5ex\hbox{$\leftrightarrow$}
\mkern-16.5mu #1}

\def\roughly#1{\raise.3ex\hbox{$#1$\kern-.75em
\lower1ex\hbox{$\sim$}}}

\def\opname#1{\mathop{\kern0pt{\rm #1}}\nolimits}

\def\dim{\opname{dim}}

\begin{document}
\quad \vskip1.375truein

\def\mq{\mathfrak{q}}
\def\mp{\mathfrak{p}}
\def\mH{\mathfrak{H}}
\def\mh{\mathfrak{h}}
\def\ma{\mathfrak{a}}
\def\ms{\mathfrak{s}}
\def\mm{\mathfrak{m}}
\def\mn{\mathfrak{n}}

\def\Hoch{{\tt Hoch}}
\def\mt{\mathfrak{t}}
\def\ml{\mathfrak{l}}
\def\mT{\mathfrak{T}}
\def\mL{\mathfrak{L}}
\def\mg{\mathfrak{g}}
\def\md{\mathfrak{d}}

\title[$C^0$-coerciveness of Moser's problem]{
$C^0$-coerciveness of Moser's problem and smoothing area
preserving homeomorphisms}

\author{Yong-Geun Oh}
\thanks{Partially supported by the NSF grant \# DMS 0503954 and a grant of the
2000 Korean Young Scientist Prize}
\date{Revision, Nov 27, 2006}

\address{
Department of Mathematics, University of Wisconsin, Madison, WI
53706 \& Korea Institute for Advanced Study, 207-43
Cheongryangni-dong Dongdaemun-gu, Seoul 130-012, KOREA,
oh@math.wisc.edu}

\begin{abstract}
In this paper, we establish the $C^0$-coerciveness
of Moser's problem of mapping one smooth volume form to another
in terms of the weak topology of measures associated to the volume forms.
The proof relies on our analysis of Dacorogna-Moser's solution to
Moser's problem of mapping one volume form to the other with
the same total mass. As an application, we give a proof of smoothing result of
area preserving homeomorphisms and its parametric version in two dimension, (or
more generally in any dimension in which the smoothing theorem
of homeomorphisms is possible, e.g., in dimension 3 but not
necessarily in dimension 4). This in turn
results in coincidence of the area-preserving homeomorphism
group and the symplectic homeomorphism group in two dimension.
\end{abstract}


\maketitle

\tableofcontents

\section{Introduction and the main theorems}
\label{sec:intro}

Consider a symplectic manifold $(X,\omega)$ and denote by
$Diff(X)$ the group of smooth diffeomorphisms of $X$.
Eliashberg's celebrated $C^0$ rigidity theorem \cite{eliash:rigid},
\cite{gromov:pseudo} in
symplectic geometry states that the subgroup $Symp(X,\omega)$
of $Diff(X)$ consisting of symplectic diffeomorphisms, i.e.,
those $\eta$ satisfying $\eta^*\omega = \omega$  is $C^0$ closed
in $Diff(X)$. More precisely, we equip the group $Homeo(X)$
of homeomorphisms with the metric
defined as
$$
\bar d(h,k) = \max_{x \in X}(d(h(x),k(x)) + d(h^{-1}(x),k^{-1}(x)))
$$
where $d$ is a distance of any given Riemannian metric. With
this metric, $Homeo(X)$ becomes a topological group which is
a complete metric space. We consider the induced topology
on $Diff(X) \subset Homeo(X)$. Eliashberg's rigidity theorem
then can be phrased as $Symp(X,\omega)$ is a \emph{closed} topological
subgroup of $Diff(X)$ with respect to this induced topology.
Motivated by this rigidity theorem, we defined
$$
Sympeo(X,\omega): = \overline{Symp}(X,\omega)
$$
where $\overline{Symp}(X,\omega)$ is the closure of $Symp(X,\omega)$
in $Homeo(X)$,  and called this group
the group of \emph{symplectic homeomorphisms} \cite{oh:hameo1}.
With this definition, the rigidity theorem can be succinctly written as
$$
Sympeo(X,\omega) \cap Diff(X) = Symp(X,\omega).
$$
Then  in the same paper \cite{oh:hameo1}, we introduced
the notion of \emph{Hamiltonian homeomorphisms} and denote the
set thereof by $Hameo(X,\omega)$. This is the $C^0$ counterpart of
the group $Ham(X,\omega)$ of Hamiltonian diffeomorphisms. We also proved
that $Hameo(X,\omega)$
forms a path-connected \emph{normal subgroup} of $Sympeo_0(X,\omega)$,
and conjectured that $Hameo(X,\omega)$ is a \emph{proper} subgroup of
$Sympeo_0(X,\omega)$. We refer readers to \cite{oh:hameo1} for further discussions
on the structure of the Hamiltonian homeomorphism group.

In two dimensional compact surface $(\Sigma,\Omega)$ with an
area form $\Omega$, we denote by $Homeo^\Omega(\Sigma)$ the group of
$\Omega$-area preserving homeomorphisms on $\Sigma$.
It easily follows from the definition that $Sympeo(\Sigma,\Omega)$
is the subgroup of $Homeo^\Omega(\Sigma)$ that consists of area preserving
homeomorphisms approximable by area preserving (smooth) diffeomorphisms.

The main motivation of the present paper
is to prove the following result conjectured in \cite{oh:hameo1}.
\medskip

\noindent{\bf Theorem I.} \emph{For a two dimensional surface $(\Sigma,\omega)$,
we also write $\omega =\Omega$ as an area form. Then we
have
$$
Sympeo(\Sigma,\omega) = Homeo^\Omega(\Sigma), \quad
Sympeo_0(\Sigma,\omega) = Homeo^\Omega_0(\Sigma).
$$
Here we denote by $G_0$ the identity component of any topological group
$G$.}
\medskip

Theorem I and normality of $Hameo(D^2,\del D^2)$ in
$Sympeo(D^2,\del D^2)$ and path-connectedness of
$Hameo(S^2,\Omega)$ proven in \cite{oh:hameo1} are the bases
on the conjecture on the structure of $Homeo^\Omega(D^2,\del D^2)$ made
in \cite{oh:hameo1}, \cite{oh:icm},
which reads that \emph{$Homeo^\Omega(D^2,\del D^2)$ is not a simple group}.

In more concrete terms, this theorem can be rephrased as the following
smoothing result of area preserving homeomorphisms which
is one belonging solely to the realm of area preserving dynamical system.
This smoothing result seems to have been a folklore among the experts in
the area but we could not locate any proper reference containing
its proof in the literature.
\medskip

\noindent{\bf Theorem $\mbox{\bf I}'$.} \emph{Let $\Sigma$ be a compact surface without
boundary and $\Omega$ be an area form. Denote by $\mu_\Omega$ the
Borel measure induced by the integration of $\Omega$. Then,
\begin{enumerate}
\item any area preserving homeomorphism $h$ can be $C^0$ approximated by
an area preserving diffeomorphism
\item
any isotopy $h = \{h_t\}_{0 \leq t \leq 1}$ with $h_0 = id$
of area preserving homeomorphisms can be $C^0$ approximated by
a smooth isotopy of area preserving diffeomorphisms.
\end{enumerate}}
\medskip
As our proof will show, Theorem $\mbox{\rm I}'$ holds for any
Borel measure induced by a volume form (or by a volume density if
not orientable) on general compact manifolds in general dimension,
\emph{as long as approximation of any homeomorphism on a manifold
$X$ by a diffeomorphism is possible}, for example in dimension 2 and 3
\cite{munkres:smoothing} but possibly not in dimension 4
\cite{donald}. It seems to be an interesting open question to ask
whether the measure preserving property helps one to approximate a
homeomorphism by a diffeomorphism and so to prove Theorem
$\mbox{\rm I}'$ in complete generality in high dimensions.

To highlight the main point of the present paper,
we outline our proof of (1) here. Denote by $M[\Sigma,\Omega]=
Homeo^\Omega(\Sigma)$ the topological group of measure preserving homeomorphisms on
$(\Sigma,\mu_\Omega)$ equipped with the topology induced by the metric $\bar d$
defined above. We call this topology the $C^0$
topology of $Homeo^\Omega(\Sigma)$. We will also denote by $d_{C^0}$ the usual
$C^0$ metric given by
$$
d_{C^0}(h,k) = \max_{x \in X} d(h(x),k(x)).
$$

Let $h \in M[\Sigma,\Omega]$ and $\e > 0$ be given.
By the well-known smoothing theorem
(see the proof of Theorem 6.3 \cite{munkres:smoothing}, for example)
for $\dim X =2$, we can choose a diffeomorphism $\psi_1$ such that
\be\label{eq:dhpsi3e}
\bar d(h,\psi_1) \leq \frac{\e}3.
\ee
This diffeomorphism $\psi_1$ however may {\it not} necessarily be
area preserving. We therefore modify $\psi_1$ into an area
preserving diffeomorphism by a $C^0$ small perturbation.

Here we would like to emphasize
that \emph{the two forms $\psi_1^*\Omega$ and $\Omega$ are not necessarily
$C^0$ close}. More precisely, we have
$$
\psi_1^*\Omega = f \Omega, \quad f > 0
$$
where $f = \det  d\psi_1$ with $d\psi_1$ being the derivative of $\psi_1$.
Since we do not have any control on $d\psi_1$ in the $C^0$
convergence, the modulus $|f-1|$ is not necessarily small. We denote
$$
|g| = \max_{x \in \Sigma}|g(x)|
$$
for a function $g$ in general. However it is not difficult to see that
(\ref{eq:dhpsi3e}) also implies that
the measures associated to $\Omega$ and $\psi_1^*\Omega$
can be made arbitrarily close in the {\it weak topology}
of measures if one chooses $\e$ sufficiently small. (See Proposition \ref{inMMX}.)

It is well-known that the set $\MM(X)$ of finite measures on a compact
metric space $X$ is a metric space such that the subset $\MM^m(X)$
of measures whose total mass is less than equal to $m \in \R_+$
is compact.
(See \cite{gromov:metric} for example.)
We denote by $d_\MM$ a corresponding metric on $\MM(X)$.
Now we will derive the proof of Theorem $\mbox{\rm I}'$
from the following theorem concerning coerciveness of the $C^0$ distance
with respect to the weak topology of measures.
\emph{This theorem holds in arbitrary dimension}. We assume $X$ is
orientable for the simplicity. Non-orientable case will be
the same if we replace the volume form by the density.
We denote by $\mu_\sigma$ the measure induced by the volume form
$\sigma$ in general.

The main result of the present paper is then the following
$C^0$-coerciveness of such diffeomorphisms $\psi_2$ in terms of
the distance $d_\MM(\mu_{\psi_1^*\Omega},\mu_\Omega)$ or in terms of
the weak topology of measures.
\medskip

\noindent{\bf Theorem II.} \emph{Let $\sigma$ and $\tau$ be two volume
forms $\sigma = f \tau$ on $X$ with $f$ satisfying $f > 0$.
Let $\lambda > 0$ be the constant
$$
\lambda = \int_X \sigma \Big\slash \int_X \tau.
$$
Then there exists a diffeomorphism $\psi_2: X \to X$ such that
$$
\psi_2^*\sigma = \lambda \tau.
$$
Furthermore, we have
\be\label{eq:dMMX}
\bar d(\psi_2,id) \to 0 \quad \mbox{ as $d_\MM(\mu_\sigma,\mu_\tau) \to 0$.}
\ee
Moreover its parametric version also holds : For any isotopy of forms
$t \in [0,1] \mapsto f_t\tau$ where $t \mapsto \mu_{(f_t\tau)}$ defines a
continuous path in $\MM(X)$, there exists an isotopy $t \in [0,1]
\to \psi_{2,t}$ of diffeomorphisms satisfying $\psi_{2,t}^*\sigma = \lambda_t \tau$
that is continuous in the compact open topology.}
\medskip

In fact, our proof of the parameterized version of Theorem II
provides canonical local slices of the action of $Homeo(X)$ on $\MM(X)$
$$
\Psi_{\mu_0} : U_{\mu_0} \cap \MM(X;\Omega) \to Homeo(X)
$$
around $\mu_0 = \mu_{g\Omega}$ with $g$ continuous,
where $U_{\mu_0} \subset \MM(X)$ is an open neighborhood of
in $\MM(X)$, and $\MM(X;\Omega)$ is the space of measures
that are absolutely continuous with respect to $\mu_\Omega$.
We will elaborate this generalization elsewhere.

Once we have Theorem II, we apply the theorem to the forms
$$
\sigma = \Omega, \quad \tau = (\psi_1^{-1})^*\Omega, \quad
\mbox{with }\, \lambda = 1
$$
and construct $\psi_2$ such that
$$
(\psi_2)^*\Omega = (\psi_1^{-1})^*\Omega \quad \mbox{and }\,
\bar d(\psi_2,id) \leq \frac{\e}{3}
$$
by letting $d_\MM(\mu_\Omega,\mu_{(\psi_1^{-1})^*\Omega})$ as small as we want.
The last can be achieved if we choose $\psi_1$ sufficiently $C^0$ close to
the area preserving homeomorphism $h$.
Then we prove that the composition $\phi:= \psi_2 \circ \psi_1: X \to X$
is an area preserving diffeomorphism with the estimate
$$
\bar d(\phi,id) \leq \e
$$
for any given $\e > 0$. A simple examination of the proof will also
give rise to the proof of its parametric version. This
will then finish the proof of Theorem II and so Theorem $\mbox{\rm I}'$.

Theorem II \emph{without (\ref{eq:dMMX})} is a result proven by
Moser \cite{moser}. And the $C^{k+1,\alpha}$ estimate for $k \geq
0, \, 0 < \alpha < 1$ that is the H\"older analog to
(\ref{eq:dMMX}) was also proven by Dacorogna and Moser \cite{DM}.
The main point of Theorem II is the $C^0$-coerciveness with
respect to the weak topology of measures which is a crucial
ingredient in our proof of the smoothing theorem, Theorem $\mbox{\rm I}'$.
We prove this coerciveness by analyzing the $C^0$-behavior
of Dacorogna-Moser's solution on the cube obtained by
their `elementary approach' with respect to the weak
topology of measures.

For this purpose, we first have to turn Dacorogna and
Moser's original {\it one-dimensional}
scheme into an \emph{$n$-dimensional} scheme which optimally reflects its
$n$-dimensional measure theoretic behavior, and to use sufficiently
small cubes whose size depends only on the given reference
volume form on $X$. (See section \ref{sec:cube}, especially Remark
\ref{Qn-X}, and the proof of Proposition \ref{|Psi0|}.) Furthermore we like to mention that,
interestingly enough, \emph{open mapping theorem}
plays an essential role in our derivation of $C^0$-coercive estimates of
Dacorogna-Moser's solution with respect to the weak topology of
measures. (See section \ref{sec:linearization}.)

Now we mention some related results in the literature. In their seminal paper, among
other things, Oxtoby and Ulam \cite{oxtoby} proved an approximation
of measure preserving homeomorphisms by \emph{almost everywhere}
differentiable measure preserving homeomorphisms. Our proof relies
on a smoothing result of general homeomorphisms for $n=2$. This
result can be extracted from \cite{munkres:S2},
\cite{munkres:smoothing} and the references therein, for example. We
refer readers to the proof of Theorem 6.3 \cite{munkres:smoothing}
in particular.
The result Theorem $\mbox{\rm I}'$ itself seems to have been a
folklore among the experts. However we have not been able to locate
a proper reference containing its proof (or its statement) in the literature. The main
theorem, Theorem II, has its own separate interest in its possible
relation to the study of generalized flows of incompressible perfect
fluids and to the problem of optimal transport (See
\cite{brenier:jams}, \cite{shnirel:gafa}, \cite{villani} for
example.)

Organization of the contents is in order. Section \ref{sec:weak}
summarizes the basic facts on the weak topology of measures relevant
to the proofs of Theorem $\mbox{\rm I}'$ and II. Section
\ref{sec:cube} recalls and enhances Moser's reduction procedure \cite{moser}
of the problem to one on the cube establishing continuity of the procedure
in the weak topology of measures. Section \ref{sec:DM} reviews
Dacorogna-Moser's elementary approach closely and provides
a reformulation of their scheme so that we can analyze its
dependence on the weak topology of measures.
After then, section \ref{sec:linearization}-\ref{sec:theoremII}
contain the proof of Theorem II.
The proof of Theorem $\mbox{\rm I}'$ will be finished in section
\ref{sec:TheoremI}.  Finally in section \ref{sec:higher}, we
prove the key a priori estimates for the study of $C^0$-coercive
estimates of Dacorogona-Moser's solution. This last section contains the most
technical estimates of the paper, whose validity, however, is
motivated by Taylor's remainder theorem.

We like to thank A. Fathi, J. Franks and J. Mather for a useful
communication during our preparation of the paper
\cite{oh:hameo1}. In reply to our question, they convinced us that
the smoothing result stated in Theorem $\mbox{\rm I}'$ should hold.
We also greatly thank the unknown referee for pointing out some serious
flaw in the previous version of the paper.
\bigskip

\noindent {\bf Notations.}
\begin{enumerate}
\item $Q = Q^n = [0,1]^n$, the unit cube in $\R^n$
\item $Q^n(K) = [0,K]^n$, the cube with its size length $K > 0$
\item For a positive number $\eta$ with $0< \eta < 1$, we denote
\beastar
Q^n(1+\eta) & = & \{x \in \R^n \mid -\eta \leq x_j \leq 1 +\eta,
\, j = 1, \cdots, n\}\\
Q^n(1-\eta) & = & \{x \in \R^n \mid \eta \leq x_j \leq 1 -\eta,
\, j = 1, \cdots, n\}
\eeastar
\item For a vector $a \in Q^n$ and $b \in [-\eta,\eta]^n$,
we denote
$$
x^{n-1}_{a;j} = (x^{j-1}, a_j, \widetilde x_j) \quad \mbox{for $j = 1,
\cdots ,n$}
$$
where we denote $\widetilde x_j: = (x^{j+1},\cdots, x^n)$.
\item $R^n_{au;j}$ : see (\ref{eq:Rnauj}).
\item $Q^n_{a;j}$ : see
(\ref{eq:Qnaj}).
\item $R^{n-1}_{au;j,k}$ : see
(\ref{eq:Rn-1ajk}).
\item $Q^{n}_{a;j,k}$ : see
(\ref{eq:Qnajk}).
\item $C^{\vec 1}(Q^n,\R^n)$ : see Definition \ref{Cvec1}.
\end{enumerate}

\section{Weak topology of $\MM(X)$}
\label{sec:weak}

In this section, we briefly review the \emph{weak topology}
of the space of finite measures on a compact metric space $X$ following
the exposition from section $3\frac{1}{2}.9$ \cite{gromov:metric}.

\begin{defn}[\bf Weak topology] A sequence of finite measures $\mu_i$ is said
to converge to $\mu$ if $\mu_i(f) - \mu(f)$ for every
bounded, nonnegative, continuous function $f$ on $X$, where $\mu(f)$
stands for $\int_X f d\mu$. We denote by $\MM(X)$ the set of
finite measures equipped with this topology.
\end{defn}

It turns out the weak topology is induced by a metric. One such
metric can be defined by
\be\label{eq:Lidb}
\operatorname{Lid}_b(\mu,\mu'): = \sup_{f}
|\mu(f) - \mu'(f)|
\ee
for $b>0$, where $f$ runs over all \emph{1-Lipschitz} functions
$f: X \to [0,b]$. These define true metrics on $\MM(X)$ and
they are mutually bi-Lipschitz equivalent. The metrics are also
complete and if $X$ is compact, then the subset of $\MM(X)$
$$
\MM_m = \{ \mu \in \MM(X) \mid \mu(X) \leq m\}
$$
is compact for each fixed $m \in \R_+$. We denote
$$
d_{\MM} = \operatorname{Lid}_1.
$$
There is a natural map
\be\label{eq:h*mu}
Homeo(X) \times \MM(X) \to \MM(X) ; (h,\mu) \mapsto h_*\mu
\ee
which is continuous (see Proposition 1.5, \cite{fathi} for example).

Next we consider the Borel measures induced by volume forms.
Let $\Omega$ be a volume form on a compact manifold $X$
satisfying $|\Omega|:= \int_X \Omega < \infty$, and denote by
$\mu_\Omega$ the measure induced by integrating the form $\Omega$.
Denoting by $\Omega^n(X)$ the space of volume forms, there is
a natural action of $Diff(X)$
\be\label{eq:psi*Omega}
Diff(X) \times \Omega^n(X) \to \Omega^n(X) ; \quad (\psi,\Omega) \mapsto
\psi^*\Omega
\ee
which is continuous in $C^\infty$ topology. It also induces a map
$$
Diff(X) \times \Omega^n(X) \to \MM(X) ; \quad (\psi,\Omega) \mapsto
\psi_*(\mu_\Omega) = \mu_{(\psi^*\Omega)}.
$$
The following proposition will play an essential role in our proof.

\begin{prop}\label{inMMX} Let $\psi \in Diff(X)$ and $h \in Homeo^\Omega(X)$
and $\psi \to h$ in $C^0$ topology. Then we have
$$
\mu_{(\psi^*\Omega)} \to \mu_\Omega \quad \mbox{ in } \MM(X).
$$
The convergence is uniform over any given compact family of $h$s.
\end{prop}
\begin{proof} It follows that $\mu_{(\psi^*\Omega)} = \psi_*(\mu_\Omega)$.
Since $\psi \to h$ in $C^0$, continuity of (\ref{eq:h*mu}) implies
$$
\psi_*(\mu_\Omega) \to h_*(\mu_\Omega)
$$
in $\MM(X)$. On the other hand we have $h_*(\mu_\Omega) = \mu_\Omega$
by the hypothesis $h \in Homeo^\Omega(X)$.
This finishes the proof of the first statement. The second statement is
an immediate consequence of the compactness assumption of the family.
\end{proof}

Now denote
$$
C^0(X,\R_+) = \{ f\in C^0(X,\R) \mid f > 0 \}
$$
and consider a volume form $\Omega$. $\Omega$ induces a natural embedding
$$
\iota_\Omega : C^0(X,\R_+) \hookrightarrow \MM(X)
$$
defined by
$$
\iota_\Omega(f): = \mu_{(f\Omega)}.
$$
This is a Lipschitz map which satisfies
\be\label{eq:lip}
d_\MM(\mu_{(f\Omega)},\mu_{(f'\Omega)}) \leq |\Omega|\cdot |f-f'|.
\ee

\section{Reduction of Theorem II to the cube}
\label{sec:cube}

In this section, we reduce the proof of Theorem II to the case of
the cube $Q = [0,1]^n \subset \R$. This reduction will be based on
a refinement of Lemma 1 \cite{moser}, Proposition \ref{local} below.
The main refinements lie in the statements (1) and (3) thereof.

Let $\{U_0,\cdots, U_m\}$ be an open covering of $X$ each element
of which can be mapped one to one onto the unit cube $Q = (0,1)^n$.

\begin{prop}[Compare with Lemma 1 \cite{moser}]\label{local}
Let $X$ be a compact manifold without boundary and let $\Omega$ and
$f \Omega$ be a volume form and $f$ a positive function satisfying
$\int \Omega = \int f \Omega$.
Then there exists decomposition of $g= f-1$
$$
g = \sum_{j=0}^m g_j
$$
where $g_j$ has support in $U_j$, and satisfies the following properties :
\begin{enumerate}
\item For all $k =0, \cdots, m$, we have
\be\label{eq:fk>0}
f_k:= 1 + \sum_{i=0}^k g_i > 0
\ee
and in particular $f_k\Omega$ defines a natural measure
$\mu_{(f_k\Omega)}$ by integrating the form $f_k\Omega$.
\item For all $k = 0, \cdots, m$,
\be\label{eq:intgk=0}
\int g_k \Omega = 0 \quad\mbox{or equivalently } \,
\int \Omega = \int f_k \Omega.
\ee
\item We have
\be\label{eq:dMMfkf}
d_\MM(\mu_{(f_k\Omega)},\mu_{(f\Omega)}) \leq C_1
\ee
where $C_1 = C_1(d_\MM(\mu_\Omega, \mu_{(f\Omega)})$ is a constant depending
on $d_\MM(\mu_\Omega, \mu_{(f\Omega)})$ and the covering only
and satisfying $C_1  \to 0$ as
$d_\MM(\mu_\Omega, \mu_{(f\Omega)}) \to 0$.
\item If $g \in C^k$, so is $g_j \in C^k, \, k\geq 0$.
\end{enumerate}
\end{prop}
\begin{proof} Following the proof of Lemma 1 \cite{moser}, we choose
a partition of unity $\phi_j \geq 0$ subordinate to the covering
$U_0,\cdots, U_m$. We order the elements $U_j$ so that for every $k = 1, \cdots, m$
the patch $U_k$ intersects $\cup_{j < k}U_j$. We denote by $\rho(k)$
any integer with $\rho(k) < k$ such that $U_k \cap U_{\rho(k)} \neq
\emptyset$. Then define the matrix $\alpha = (\alpha_{jk})$
with $0 \leq j \leq m$ and $1 \leq k \leq m$ by
$$
\alpha_{jk} = \begin{cases}  1 &\quad \mbox{for $j=k$} \\
 -1 &\quad  \mbox{for $j = \rho(k)$}, \\
 0 &\quad  \mbox{otherwise}.
\end{cases}
$$
This matrix satisfies $\sum_{j=0}^m\alpha_{jk} = 0$.

We now fix functions $\eta_k, \, k = 1, \cdots, m$ such that
\be\label{eq:intetak}
\int \eta_k \Omega = 1.
\ee
We can choose them so that
$$
|\eta_k| \leq C_2
$$
where $C_2$ depends only on the covering and $\Omega$. We will represent
$g_j$ in the form
$$
g_j = g \phi_j - \sum_{k=1}^m \lambda_k\alpha_{jk}\eta_k.
$$
Then Moser \cite{moser} showed that $g_j$ is add up to $g$ and has
support in $U_j$. To prove (\ref{eq:intgk=0}), we consider
the linear equation
\be\label{eq:lambdak}
\sum_{k=1}^m \lambda_k\alpha_{jk} = \int_X (f-1) \phi_j\Omega
\ee
for $j = 0, \cdots, m$, which has $m$ unknowns and $m+1$ equations.
However, on account of (\ref{eq:intetak}) and the equation
$$
\int (f-1) \Omega = \sum_{j=0}^m \int_X (f-1)\phi_j\Omega = 0
$$
the first equation ($j=0$) of (\ref{eq:lambdak}) is redundant. Therefore the solution
space of (\ref{eq:lambdak}) is a nonempty affine subspace of $\R^m$.
So far our proof has been a duplication of Moser's \cite{moser}.

The new statements in this proposition that were not
considered in \cite{moser} or \cite{DM}
are (\ref{eq:fk>0}) and (\ref{eq:dMMfkf}).
To establish these statements, we
need to analyze the solution space of (\ref{eq:lambdak})
more closely than \cite{moser} or \cite{DM} do in terms of
the weak topology of measures. First we note that we have
$$
f_j = 1 + \sum_{i=0}^j \left(\phi_i (f-1) - \sum_{k=1}^m
\lambda_k\alpha_{ik}\eta_k \right).
$$
We set $f_{-1} = 1$ and note $f_m = f$. Thanks to (\ref{eq:lip}),
to prove (\ref{eq:fk>0}) and (\ref{eq:dMMfkf}),
it will be enough to make the norms $|\lambda_k|$ all sufficiently
small. To be more precise, we rewrite (\ref{eq:lambdak}) into
$$
\sum_{k=1}^m \lambda_k\alpha_{jk} =
\int_X \phi_j d\mu_{(f\Omega)} - \int_X \phi_jd\mu_\Omega
$$
for $j = 0, \cdots, m$. Recalling the definition of $d_\MM = \operatorname{Lid}_1$,
note that the right hand side is bounded by
$$
\Big|\int_X \phi_j d\mu_{(f\Omega)} - \int_X \phi_jd\mu_\Omega\Big|
\leq d_\MM(\mu_{(f\Omega)},\mu_\Omega)
$$
since $\phi_j$ is a function satisfying $0 \leq \phi_j \leq 1$.
We like to alert the readers that \emph{the distance in the right side of
this inequality is in terms of the weak topology of measure}.

A simple linear algebra then concludes that there exist solutions
$\lambda_k$ of (\ref{eq:lambdak}) such that
\be\label{eq:|lambdak|}
|\lambda_k| \leq C_3,\quad k = 1, \cdots, m
\ee
where $C_3=C_3(d_\MM(\mu_{(f\Omega)},\mu_\Omega))$ is a constant
depending only on $d_\MM(\mu_{(f\Omega)},\mu_\Omega)$
such that $C_3 \to 0$ as $d_\MM(\mu_{(f\Omega)},\mu_\Omega) \to 0$ :
Note that the solution space of (\ref{eq:lambdak}) is a nonempty
affine subspace of $\R^m$, whose distance from the origin
converges to zero as $d_\MM(\mu_{(f\Omega)},\mu_\Omega) \to 0$.
To obtain such a solution $(\lambda_1, \cdots,\lambda_m)$ satisfying
(\ref{eq:|lambdak|}), one may take
the point  nearest to the origin among the points in the affine space.

Now to prove (\ref{eq:fk>0}), we consider the convex combinations
of $f$ and $1$
$$
f_j':= 1 + \left(\sum_{i=0}^j \phi_i\right)(f-1)=
\left(\sum_{i=0}^j \phi_i\right)f + \left(1- \left(\sum_{i=0}^j \phi_i\right)\right)1
$$
for $j = 0, \cdots, m$ and denote
\beastar
m_0 & = & \min_j\{\min f_j' \mid  j = 0, \cdots, m \}\nonumber\\
M_0 & = & \max_j\{\max f_j' \mid  j = 0, \cdots, m \}
\eeastar
Note that $m_0, \, M_0$ depends only on $f$ and satisfies
\be\label{eq:m0M0}
m_0 \geq \min\{\min f, 1\}, \quad M_0 \leq \max\{\max f,1\}
\ee

Recalling that $C_3 \to 0$ as $d_\MM(\mu_{(f\Omega)},\mu_\Omega) \to 0$,
we can choose $d_\MM(\mu_{f\Omega},\mu_\Omega)$ so small that we have
$$
C_3 < \frac{\min\{\min f,1\}}{m(m+1) C_2}.
$$
Then we derive
\be\label{eq:minminf1}
\max\sum_{k=1}^m |\lambda_k||\eta_k|\leq m C_2C_3 < \frac{\min\{\min f,1\}}{m+1}
\ee
from (\ref{eq:|lambdak|}). Therefore we  have
\beastar
f_j & = & f_j' - \sum_{i=0}^j\sum_{k=1}^m \lambda_k\alpha_{ik}\eta_k
\geq f_j' - \sum_{i=0}^j\sum_{k=1}^m |\lambda_k||\eta_k|\\
& \geq & f_j' - (m+1) \sum_{k=1}^m |\lambda_k||\eta_k|
> f_j' - \min\{\min f,1\}\\
& \geq & M_0 - \min\{\min f,1\} \geq 0
\eeastar
which proves (\ref{eq:fk>0}).

Finally we consider $d_\MM(\mu_{(f_k\Omega)},\mu_{(f\Omega)})$ for
the proof of (\ref{eq:dMMfkf}).
Since $f_k'$ is a convex combination of $f$ and $1$, we have
\be\label{eq:dMMfk'f}
d_\MM(\mu_{(f_k'\Omega)},\mu_{(f\Omega)}) \leq
d_\MM(\mu_\Omega, \mu_{(f\Omega)})
\ee
for any $k = 1, \cdots, m$.  We have
$$
d_\MM(\mu_{(f_k'\Omega)},\mu_{(f_k\Omega)}) \leq |\Omega|\cdot|f_k' - f_k|
$$
from (\ref{eq:lip}) and
$$
|f_k'-f_k| = \Big|\sum_{k=1}^m \lambda_k\alpha_{jk}\eta_k\Big|
\leq \sum_{k=1}^m |\lambda_k|\cdot |\eta_k|.
$$
On the other hand, we can make
$\max (\sum_{k=1}^m |\lambda_k||\eta_k|)$
as small as we want by choosing $\lambda_k$ small which
in turn can be achieved by (\ref{eq:|lambdak|}) if we make
$d_\MM(\mu_{(f\Omega)},\mu_\Omega)$ sufficiently small.
This implies that we can also make
$d_\MM(\mu_{f_k'\Omega},\mu_{f_k\Omega})$ as small as we want if
we make $d_\MM(\mu_{(f\Omega)},\mu_\Omega)$ sufficiently small.
We note the triangle inequality
$$
d_\MM(\mu_{(f_k\Omega)},\mu_{(f\Omega)}) \leq d_\MM(\mu_{(f_k\Omega)},\mu_{(f_k'\Omega)})
+ d_\MM(\mu_{(f_k'\Omega)},\mu_{(f\Omega)}).
$$The last statement of the proposition is obvious from the construction of $g_j$'s.
This finishes the proof.
\end{proof}

With Proposition \ref{local} in our hand, Theorem II will be derived
from the following proposition. \emph{Except the coerciveness
(\ref{eq:MMQ}), this is precisely Lemma 2 \cite{moser} or
Proposition 8 \cite{DM}}. However our diffeomorphism
may not necessarily the same as the one constructed in \cite{DM}.
In fact, our construction will provide
\emph{continuous local slices} under the action
$$
Homeo(X) \times \MM(X) \to \MM(X)
$$
over a certain dense subset of $\MM(X)$. We will elaborate this
generalization elsewhere.

\begin{thm}\label{c0est} Let $Q$ be the square $[0,1]^n$.
Consider two volume forms
$$
\tau = f(x) dx,\quad \sigma = g(x) dx
$$
where $g, \, f$ are positive continuous functions for which
$g-f$ has support in $\text{Int } Q$. Denote by $m_f, \, m_g$
the associated measures. If
\be\label{eq:fg}
\int_Q f\, dx = \int_Q g\, dx
\ee
then there exists a diffeomorphisms $\psi: Q \to Q$ such that
\be\label{eq:=hx}
g(\psi(x)) \operatorname{det} \nabla\psi(x) = f(x)
\ee
such that $\psi(x) = x$ near the boundary of $Q$.
Furthermore $\psi$ satisfies the following additional properties :
\begin{enumerate}
\item We can make $\bar d(\psi,id)$ as small as we want by letting
$d_\MM(m_f,m_g) \to 0$, or
\be\label{eq:MMQ}
\bar d(\psi,id) \to 0 \quad \mbox{ as $d_\MM(m_f,m_g) \to 0$.}
\ee
And the parametric version in the sense as stated in Theorem II also holds.
\item
Let $\operatorname{supp}(\psi) = \overline{\{x \in Q \mid \psi(x) \neq x\}}$.
Let $R^n \subset (0,1)^n$ be any closed cube such that
\be\label{eq:Rn}
R^n \subset \operatorname{supp}(f-g).
\ee
Then we have
\be\label{eq:supp}
\operatorname{supp}(\psi_2) \subset R^n
\ee
\end{enumerate}
\end{thm}

\begin{rem} \label{Qn-X}
\begin{enumerate}
\item
Obviously, we can further decompose the cube $[0,1]^n$ or
use cubes of the smaller size in Proposition \ref{local}, and get the
same kind of statement for the smaller cubes. Later in our
estimates, we will need to choose a cube $Q^n(K)$ of its side length
$K > 0$ such that $K$ is sufficiently small and depends essentially
on the given fixed $g$. In fact, we can choose $K$ of the form
$K = 2^{-N_0}$ with
$$
2^{-N_0} < \frac{1}{8C_4 (1+L_g)}
$$
where $C_4 = \max\{8 \max g,4\}$ and $L_g$ is the modulus of continuity of
$g$. See the paragraph around (\ref{eq:C4}) for more discussion on this.
However to make our exposition better comparable to that of
\cite{DM}, we will carry our discussion on the unit cube and just
indicate the needed changes in the paragraph around (\ref{eq:C4}).
\item We also note that the above reduction procedure to the cube
shows that the distance $d_{\MM}(m_f,m_g)$ for the measures
$m_f, \, m_g$ on $Q^n(K)$ converges to zero uniformly as
$d_{\MM}(\mu_\tau,\mu_\sigma) \to 0$ for the originally given
measures $\mu_\tau,\, \mu_\sigma$ on $X$.
\item The inequality (\ref{eq:m0M0}) shows that the above reduction
procedure essentially does not decrease the lower bound $\min g$
and not increase the upper bound $\max g$ on the cube from that
of the originally given $g$ on $X$.
\emph{We would like to warn the readers that $f$ in the proof of
Proposition \ref{local} plays the role of $g$ in Theorem \ref{c0est}
and henceforth}.
\end{enumerate}
\end{rem}

The next three sections will be occupied by the proof of statement (1)
of this theorem.

\medskip

\section{Scheme of construction on the cube}
\label{sec:DM}

In this section, we first recall Moser's or Dacorogna and Moser's
`elementary approach' from \cite {moser}, \cite{DM} to solving
(\ref{eq:=hx}). After then we reformulate their scheme into an
$n$-dimensional scheme so that we can study its relevance to the
weak topology of $n$-dimensional measures. Their inductive
one-dimensional approach as it is does not manifest the relationship
of their solutions with the weak topology of measures. We also
briefly mention the parametric version of Dacorogna-Moser's approach
which is used in the proof of $Sympeo_0(\Sigma,\omega) =
Homeo^\Omega_0(\Sigma)$ in Theorem I.

We denote $Q = Q^n = [0,1]^n$
and by $Q^s$ the $s$ dimensional cube for $1 \leq s \leq n$.

\subsection{Review of Dacorogna and Moser's elementary approach}
\label{subsec:review}

According to Dacorogna and Moser \cite{DM}, under the assumption
as in Theorem 3.2 on $f$ and $g$, the map
$\psi:Q \to Q$ is constructed as the successive composition
$$
\psi = \varphi_n \circ \varphi_{n-1}\circ \cdots \circ
\varphi_2 \circ \varphi_1
$$
by defining $g_n = g$ and for $s = 2,3,\cdots,n$ and requiring
\be\label{eq:gs-1gs}
\int_E g_{s-1}(x)\, dx = \int_{\varphi_s(E)}g_s(x)\, dx
\ee
for every open set $E \subset Q^n$ and
\be\label{eq:g1f}
\int_0^1 g_1(x_1,x')\, dx_1 = \int_0^1 f(x_1,x')\, dx_1.
\ee
And $\varphi_1: Q \to Q$ will then have the form
$$
\varphi_1: (x_1, x_2, \cdots, x_n) \to (v(x), x_2, \cdots, x_n)
$$
where $v: Q \to Q$ is uniquely determined by the requirement
\be\label{eq:fg1v}
\int_0^a f(x_1,x')\, dx_1 = \int_0^{v(a,x')}g_1(x_1,x')\,dx_1
\ee
for every $x' = (x_2,\cdots, x_n) \in Q^{n-1}$. Since $g_1 > 0$
(\ref{eq:fg1v}) uniquely determines $v(x)$ with $v$ monotone in
$x_1$, $v =0$ for $x_1 =0$ and $v = x_1$ for $x'$ near
$\del Q^{n-1}$. Finally (\ref{eq:g1f}) makes $v(1,x') = 1$
for all $x' \in Q^{n-1}$.
It follows that $g_s \in C^k, \, k\geq 1$, (\ref{eq:gs-1gs})
is equivalent to
\be\label{eq:gs-1=}
g_{s-1}(x) = g_s(\varphi_s(x))\det \nabla \varphi_s(x).
\ee

Then they construct
$\varphi_n, \cdots, \varphi_2$ (and $g_{n-1}, \cdots, g_1$)
inductively in such a way that
\be\label{eq:gfQs}
\int_{Q^s} g_s(x^s,x')\, dx^s = \int_{Q^s} f(x^s,x')\, dx^s
\ee
where $x^s = (x_1, \cdots, x_s)$ and $x'=(x_{s+1},\cdots,x_n)$.
Assuming that $\varphi_n,\cdots,\varphi_{s+1}$ are already
constructed so that (\ref{eq:gs-1gs}) and (\ref{eq:gfQs})
hold and that they agree with the identity near the boundary,
the map $\varphi_s: Q \to Q$ is constructed as the homeomorphism of the form
\be\label{eq:varphis}
\varphi_s(x_1, \cdots, x_n) = (x^{s-1},v(x),x')
= (x_1, \cdots, x_{s-1}, v(x), x_{s+1}, \cdots, x_n)
\ee
with
\be
v(x) = x_s + \zeta(x^{s-1})u(x_s,x').
\ee
Here $\zeta$ is a cut-off function with $\text{supp }\zeta
\subset \text{Int }Q^{s-1}$ and satisfying
\be\label{eq:zeta}
\begin{cases} 0 \leq \zeta \leq 1 + \e \quad \text{in } \, Q^{s-1} \\
\int_{Q^{s-1}}\zeta(x^{s-1}) \, dx^{s-1} = 1 \\
\int_{Q^{s-1}}|\zeta(x^{s-1}) - 1| \, dx^{s-1} < \e
\end{cases}
\ee
where $\e = \e(g_s,f) > 0 $ is chosen so that
\be\label{eq:epsilon}
\e \max g_s < \min g_s, \, \frac{1}{2}\min f.
\ee
And $u:[0,1] \to [0,1]$ is a smooth function with
$$
u \equiv 0 \quad \text{near $\{0, 1\}$}.
$$
Note that in this construction the variable $x'$ enters only as a parameter
and does not play any role in finding $\varphi_s$.
Therefore we drop $x'$ in our discussion below writing
$u(x_s) = u(x_s;x')$ as in \cite{DM}.
We refer readers to (4) and (5) \cite{DM} for more details.
It follows that $\varphi_s$ is $C^0$ close to identity if and
only if the one variable function $u: [0,1] \to [0,1]$ is
$C^0$ close to the zero function. Furthermore it becomes a
differentiable homeomorphism if and only if $u$ is differentiable and
satisfies
\be\label{eq:dvdxs}
\frac{\del v}{\del x_s} = 1 + \zeta(x^{s-1}) \frac{\del u}{\del x_s}
> 0.
\ee

To solve (\ref{eq:gfQs}), Dacorogna and Moser transformed
it into the functional equation
\be\label{eq:GF}
G(x_s,u(x_s)) = F(x_s)
\ee
where
\beastar
G(a,b) & = &\int_{R_{ab}^s}g_s(x^s)\, dx^s\\
F(a) & = &\int_{Q_a^s} f(x^s) \, dx^s
\eeastar
with
\bea
Q_a^s & = & \{ x^s \in Q^s \mid  0 < x_s < a \}\label{eq:Qas}\\
R_{ab}^s & = & \{ x^s \in Q^s \mid 0 < x_s < a + \zeta(x^{s-1}) b \}
\nonumber
\eea
where $b \in  [-\eta, \eta]$ :
They obtained this equation by first setting $u(0,x') = 0$ and then
integrating the equation (\ref{eq:gfQs}) for $(s-1)$ in place of $s$, i.e.,
\be\label{eq:xs-1xs}
\int_{Q^{s-1}}(g_{s-1}(x^{s-1},x_s,s') - f(x^{s-1},x_s,x') ) \, dx^{s-1} = 0
\ee
over $0 < x_s < a$, which gives rise to
$$
\int_{Q^s_a}(g_{s-1}(x^s,x') - f(x^s,x')) \, dx^s = 0.
$$
But this is then equivalent to (\ref{eq:GF}).

We note that $G$ (resp. $F$) is differentiable, if $g$
(resp. $f$) is continuous. In fact, we have the explicit formulae
\bea
\frac{\partial G}{\partial b} & = & \int_{Q^{s-1}}
\zeta(x^{s-1})  g_s(x^{s-1},a+\zeta(x^{s-1}) b)\, dx^{s-1} \label{eq:dGdb}\\
\frac{\partial G}{\partial a} & = &\int_{Q^{s-1}}
g_s(x^{s-1},a+\zeta(x^{s-1}) b)\, dx^{s-1} \label{eq:dGda} \\
\frac{\partial F}{\partial a} & = &\int_{Q^{s-1}}
f(x^{s-1},a)\, dx^{s-1}. \label{eq:dFda}
\eea
Note that $u(0) = 0$ is the unique solution of (\ref{eq:GF}) at $x_s = 0$.
At this point, they derived existence and uniqueness
of the solution to (\ref{eq:GF}) by the intermediate value theorem.
We denote by $v = v_{DM}$ and $u = u_{DM}$ for this unique solution
and call them Darcorogna-Moser's solution, or simply as DM-solutions.

\begin{rem}
\noindent To obtain the $C^0$ convergence statement (\ref{eq:MMQ})
in Theorem \ref{c0est}, we need to control the $C^0$ distance $\bar
d(\psi_2,id)$ in the above existence proof of $\psi_2$. This $C^0$
estimate is precisely the one left untreated by Dacorogna and Moser
in \cite{DM}. However, following Moser's deformation method
\cite{moser} and the use of elliptic second order partial
differential equation, they proved an existence of a diffeomorphism
$\psi_2$ satisfying an a priori $C^{k+1,\alpha}$ estimate when $f,
\, g \in C^{k,\alpha}$ \emph{when $k \geq 1$ and $\alpha > 0$
\cite{DM}. This elliptic approach using the deformation method does
not produce the $C^0$ convergence required in (\ref{eq:MMQ})}.
\end{rem}

In fact by differentiating (\ref{eq:GF}), one obtains
\be\label{eq:dGda}
\frac{\del G}{\del a} + \frac{\del G}{\del b}
\frac{\del u}{\del x_s} = \frac{\del F}{\del a}.
\ee
From this, Dacorogna-Moser \cite{DM} derives that the solution $u$ is
differentiable. In fact, the standard boot-strap argument,
using (\ref{eq:dGda}) and the fact that the function $\frac{\del G}{\del b}$
is positive from (\ref{eq:dGdb})
proves the following a priori $C^{k,\alpha}$ estimate for $k \geq 0$
and $0 < \alpha < 1$ for the DM-solution itself. This demonstrates
that the DM-solution is as good as the one obtained by the deformation
approach used in \cite{moser}, \cite{DM} even for the higher regularity.

One main theorem we prove in the current paper is
that DM-solutions will also satisfy
the additional $C^0$-coerciveness property under the distance
$d_\MM(m_g,m_f)$ of the weak topology of measures.

For the purpose of our later study of the parametric version
of Theorem \ref{c0est}, we summarize the above discussion
on the higher regularity into the following proposition

\begin{prop} \label{Ckestimate}
Let $g$ be a given positive $C^{k,\alpha}$ function.
Suppose that the functions $f$ is also $C^{k,\alpha}$ and
denote by $|\cdots|_{k,\alpha} $ the $C^{k,\alpha}$ norm of functions.
Let $u$ be a DM-solution. Then we have
\be\label{eq:Ckestimate}
|u|_{k+1,\alpha} \leq C_{(k;g)} |f-g|_{k,\alpha}
\ee
for all $k \geq 0$, where $C_{(k;g)}$ is a constant depending only on $k$ and
$C^k$ norm of $g$.
\end{prop}

\begin{rem}
We would like to emphasize that we cannot expect that the derivative
of the solution $u$ converges to 0 as $d_\MM(m_f,m_g) \to 0$.
In fact in the above proof, we do not have any control of
$|\nabla u|$ in terms of $d_\MM(m_f,m_g)$.
\end{rem}

\subsection{Coercive reformulation}
\label{subsec:refinement}

At the end of the day, one can write Dacorogna-Moser's solution in
the form $\psi = \varphi_n \circ \varphi_{n-1} \circ \cdots \circ
\varphi_2 \circ \varphi_1$. In coordinate expression $\psi = (v_1,
v_2,\cdots, v_n):= v$, $v_j$ has the following form :

\bea\label{eq:DMpsi}
v_1(x) & = & x_1 + u_1(x_1,\widetilde x_1) \nonumber\\
v_2(x) & = & x_2 + \zeta_2(v^1(x))u_2(x_2,\widetilde x_2) \nonumber \\
& \vdots & \nonumber \\ v_{n-1}(x) & = & x_{n-1} +
\zeta_{n-1}(v^{n-2}(x))u_{n-1}(x_{n-1},x_n) \nonumber \\
v_n(x) & = & x_n + \zeta_n(v^{n-1}(x))u_n(x_n).
\eea
Here we denote
$\widetilde x_i = (x_{i+1}, \cdots,x_n)$ and $v^j = (v_1, \cdots,
v_j)$ for $j = 1, \cdots, n$. We would like to \emph{emphasize} that the argument
inside $\zeta_j$ is $v^{j-1}(x)$, not $x^{j-1}$.

We will now examine the $C^0$-behavior of DM-solutions
$\psi$ above in terms of the weak topology of measures.

We recall that $\operatorname{supp}(f-g) \subset \text{Int }Q$ and
so we can choose $\eta > 0$ so that \be\label{eq:eta}
\operatorname{supp}(f-g) \subset \{x \in Q \mid d(x,\del Q) \geq
\eta\} = Q^n(1-\eta). \ee This choice of $\eta$ depends only on
$\operatorname{supp}(f-g)$, independent of individual $f$ or $g$.
The choice of $\eta$  will be fixed for the rest of the paper.
Without loss of generality, we also assume that $f,\, g$ are indeed
defined on the bigger cube $Q^n(1 + \eta)$ where
$$
Q^n(1+\eta) = \{x \in \R^n \mid -\eta \leq x_j \leq 1 +\eta,
\, j = 1, \cdots, n\}.
$$
We now fix a family of cut-off functions $\zeta = \{ \zeta_s\}_{s =
2}^n$ with

\be\label{eq:zetas-supp} \zeta_s: Q^{s-1} \to \R \quad \mbox{with
$\operatorname{supp}\zeta_s \subset Q^{s-1}(1-\frac{\eta}{2})$} \ee
for $s = 2, \cdots, n$ such that \be\label{eq:zetas}
\begin{cases} 0 \leq \zeta_s \leq 1 + \e_0 \quad \text{in } \, Q^n \\
\int_{Q^{s-1}}\zeta_s(x^{s-1}) \, dx^{s-1} = 1 \\
\int_{Q^{s-1}}|\zeta_s(x^{s-1}) - 1| \, dx^{s-1} < \e_0
\end{cases}
\ee as in (\ref{eq:zeta}) where $\e_0 = \e_0(\zeta) > 0$
is a constant, which satisfies
\be\label{eq:e0} \e_0(\zeta) <
\min\left\{ \frac{\min\{\min f, \min g\}}{\max g}, \frac{\min
f}{2\max g} \right\}.
\ee
\emph{This constant $\e_0(\zeta)$ can be made as small as we want
independently of the given $g, \,f$}. (See Remark \ref{Qn-X} (2) and (3).)
For example, we can always choose
\be\label{eq:e0<dMM}
\e_0(\zeta) < d_\MM(m_g,m_f).
\ee
Motivated by the expression given in (\ref{eq:DMpsi}), we introduce the
following definition which will be essential for our discussion following
afterwards.
\begin{defn}\label{triangular}
We call a map $u: Q^n \to \R^n$ \emph{triangular} if its components
$u_j$ have the following triangular form :
\beastar
u_1 & = & u_1(x_1, \cdots, x_n)\\
u_2 & = & u_2(x_2,\cdots, x_n) \\
& \vdots & \\
u_{n-1} & = & u_{n-1}(x_{n-1},x_n) \\
u_n & = & u_n(x_n).
\eeastar
We denote by $C^0_{tri}(Q^n,\R^n)$ the set of triangular maps.
We define
$$
B \subset C^0_{tri}(Q^n, \R^n)
$$
the set of triangular maps satisfying $u(1,\cdots,1) = 0$.
\end{defn}

Obviously $B$ is a closed subspace of the Banach space $C^0(Q^n,
\R^n)$ and hence itself a Banach space with the $C^0$-norm
$$
|u|=|u|_{C^0} = \max_{1 \leq j \leq n}|u_j|
$$
for the vector map $u = (u_1,\cdots, u_n)$. Furthermore it follows from
this triangularity of $u$ that
the Jacobian $\nabla u$ of $u$ forms an upper triangular matrix.

Now the DM-solutions $v=v_{DM}: Q^n \to Q^n$ have the
following form

\be\label{eq:etas} v(x_1, \cdots, x_n)  = (v_1, \cdots,
v_j,\cdots, v_n) \ee where $v_j: Q \to [0,1]$ is a function of the
type

\bea v_j(x)  &  =  & v_j(x_j,\widetilde x_j)  =  x_j +
\zeta_j(v^{j-1}(x))u_j(x_j,x'),
\quad j = 2, \cdots, n \label{eq:vj} \\
v_1(x) & = & v_1(x_1, \cdots x_n) = x_1 + u_1(x_1,\cdots, x_n)
\label{eq:v1}
\eea
with $\widetilde x_j=(x_{j+1}, \cdots, x_n)$. In
other words, we can factorize $\psi_2$ into
$$
\psi_2 = \varphi_n \circ \varphi_{n-1}\circ \cdots \circ \varphi_1
$$
where each $\varphi_j$ is a smooth map of the form given in
(\ref{eq:varphis}) depending on $u$.

Then the diffeomorphism $v$ satisfies $g(v(x))\det\nabla v(x) = f(x)$
and its weak form
\be\label{eq:gfQk}
\int_{v(E)} g(y)\, dy = \int_{E} f(x)\, dx
\ee
for any measurable subset $E$. We define

\bea Q^n_{a;j} & = & \{x \in Q^n \mid 0 \leq x_i \leq 1\, \text{for }\,
1 \leq i \leq j-1, \, \nonumber \\
& {} & \hskip0.5in 0 \leq x_i \leq a_i, \, j \leq i \leq n \}
\label{eq:Qnaj}\\
R^n_{au;j} & = & v(Q^n_{a;j}) \label{eq:Rnauj}
\eea for  $j = 1, \cdots, n$.

Knowing that the DM-solution $v=\psi_2$ is a homeomorphism
(in fact a smooth diffeomorphism when $g,\, f$ are smooth),
$R^n_{au;j}$ is a closed measurable subset and so
we can define the integrals
\bea
G_j(a;u) & = & \int_{R^n_{au;j}} g(y) \, dy \label{eq:Gabj}\\
F_j(a) & = & \int_{Q^n_{a;j}} f(x) \, dx \label{eq:Faj}
\eea
and consider the vector functions
$$
G = (G_1, \cdots, G_n), \, \quad F= (F_1,\cdots, F_n)
$$
where we denote $G : = G(\cdot;u)$. Then the weak form (\ref{eq:GF})
of the equation
$$
g(v(x))\det\nabla v(x) = f(x)
$$
can be reduced to (\ref{eq:GFy})
\be \label{eq:GFy} G(a;u) = F(a), \quad a \in
Q^n.
\ee
In particular, DM-solution satisfies (\ref{eq:GFy}).

The converse also holds for differentiable maps.
\begin{lem}
If $u$ is a solution of (\ref{eq:GFy}) that is
differentiable, then it satisfies
\be\label{eq:gvx}
g(v(x))\det\nabla v(x) = f(x).
\ee
\end{lem}
\begin{proof} Since $u$ is differentiable, we can apply
the change of variables and rewrite
(\ref{eq:GFy}) as
$$
\int_{Q^n_{a;j}} g(\psi_2(x)) \nabla \psi_2(x) \, dx
= \int_{Q^n_{a;j}} f(x)\, dx
$$
for all $j = 1,\cdots, n$.
The lemma then follows by taking the partial derivatives of
these equations with respect to $a_i$ for each $i = 1, \cdots, n$.
\end{proof}

Now we consider the subset $B_{homeo} \subset B$ defined by
$$
B_{homeo} = \{ u \in B \mid \mbox{the associated map $v$ in
(\ref{eq:vj}) and (\ref{eq:v1}) is a homeomorphism} \}.
$$
Then for each element $u \in B_{homeo}$, the functions $G_j$
are defined and so we can define a map
\be
\label{eq:Psi} \Psi :  B_{homeo} \to C^0(Q^n,\R^n)
\ee
by $\Psi = (\Psi_1, \cdots, \Psi_n)$
whose components are given by
$$
\Psi_j(u) = G_j(\cdot;u) - F_j(\cdot).
$$
We remark that the equation (\ref{eq:GFy}) is equivalent to
$\Psi(u) = 0$.

The following proposition is the reason why we introduce the
notion of triangular maps and the space $B$.

\begin{prop}\label{Xi} The map $a \mapsto F(a)$ is triangular, and
so is $a \mapsto G(a;u)$ whenever $u \in B_{homeo}$.
In particular, the map $\Psi$ maps $B_{homeo}$ to $B$.
\end{prop}
\begin{proof}
Recall the definitions of $G$ and $F$ in (\ref{eq:Gabj}) and
(\ref{eq:Faj}) respectively. By the definition (\ref{eq:Qnaj}) of
$Q^n_{a;s}$, it does not depend on $a_1, \cdots, a_{s-1}$ and
hence neither does $R_{au;s}^n = v(Q_{a;s}^n)$. This immediately
implies that both $F$ and $G(\cdot;u)$ are triangular.
This finishes the proof of triangularity of $\Psi$.

We next check $\Psi(u)(1,\cdots, 1) = 0$. Since $u(1,\cdots, 1) =
0$ we have
$$
R_{(\vec 1u;n)}^n = R_{(\vec 1 \vec 0;n)}^n = Q^n.
$$
where $\vec 1 = (1,\cdots, 1)$ and $\vec 0 = (0,\cdots, 0)$.
Therefore we have
$$
\Psi(u)(1,\cdots, 1) = \int_{Q^n} g \,dy - \int_{Q^n} f \, dx
$$
which is assumed to be zero in (\ref{eq:fg}). This finishes the
proof.
\end{proof}

\section{Linearization}
\label{sec:linearization}

Now we introduce the subset $B_{diff} \subset B_{homeo}$ consisting
of smooth maps $u$ whose associated map $v$ is a
diffeomorphism. Then the restriction of $\Psi$ to
$B_{diff}$ is continuously differentiable map to $C^\infty(Q^n,\R^n)$
in the Frechet sense : Since diffeomorphism property of a map
defined on compact sets is an open property, once we know that
$B_{diff}$ is non-empty, it is an open subset of $C^\infty(Q^n,\R^n)$
and hence we can define the Frechet derivative of $\Psi$ on $B_{diff}$.

Denote by $\overline 0$ the zero function.
We now compute the Frechet derivative of $\Psi_{B_{diff}}$ at
$u= \overline 0 \in B_{diff}$ which corresponds to $v = id$.

Applying the Taylor expansion to $\Psi$ at $u = \overline 0$,
(\ref{eq:GF}) can be rewritten as
\be\label{eq:-dPsi0}
- d\Psi(\overline 0) \cdot u = \Psi(\overline 0) + N(u)
\ee
where $d\Psi$ is the Frechet derivative of
$$
\Psi: B_{diff} \to  C^\infty(Q^n,\R^n)
$$
and
$$
N(u) = \Psi(u) - \Psi(\overline 0) - d\Psi(\overline 0) \cdot u
$$
is the `higher order term'.
It follows from the definitions of $\Psi_j$ that we have
\bea
\Psi_j(\overline 0)(a) & = & G_j(a;\overline 0) - F_j(a)
\nonumber \\
& = & \int_{Q^n_{a;j}} g(y) dy - \int_{Q^n_{a;j}} f(x) dx.
\label{eq:Psia0}
\eea
Now the following provides an explicit formula for the
Frechet derivative of the map
$$
d\Psi(\overline 0): C^\infty_{tri}(Q^n,\R^n)
\to C^\infty_{tri}(Q^n,\R^n)
$$
at $u = \overline 0$.

\begin{prop} Let $X = (X_1, \cdots, X_n) \in C^\infty_{tri}(Q^n,\R^n)$. Then
\be
(d\Psi(\overline 0)\cdot X)_j(a)
= \sum_{k=j}^n \int_{Q^{n-k}_{a;j}} X_k(a_k,\widetilde x_k)
\Big(\int_{Q^{k-1}_{a;j,k}} \zeta_j(x^{k-1}) g(x^{n-1}_{a;k}) dx^{k-1}\Big)
d\widetilde x_k.
\label{eq:dPsi0}
\ee
In particular, the matrix elements
$$
(d\Psi(\overline 0))_{jk} : C^\infty(Q^{n-k},\R) \to C^\infty(Q^{n-j},\R)
$$
of the matrix operator
$$
d\Psi(\overline 0): C^\infty_{tri}(Q^n,\R^n)
\to C^\infty_{tri}(Q^n,\R^n)
$$
are given by
\be\label{eq:matrix}
\left((d\Psi(\overline 0))_{jk}(h)\right)(a) = \int_{Q^{n-k}_{a;j}}
C_{jk}(a_j,\cdots, a_k,\widetilde x_k) h(\widetilde x_k)\, d\widetilde x_k
\ee
where $C_{jk}(a_j,\cdots, a_k,\widetilde x_k)$ are smooth functions
of $(a_j,\cdots, a_k, \widetilde x_k)$ defined by
\be\label{eq:Cjk}
C_{jk}(a_j,\cdots, a_k,,\widetilde x_k) = \begin{cases}\int_{Q^{n-1}_{a;j}}\zeta_n(x^{n-1})
g(x^{n-1}, a_n) dx^{n-1} \quad &\mbox{for }\, k = n\\
\int_{Q^{k-1}_{a;j}}\zeta_k(x^{k-1})
g(x^{k-1}, a_k, \widetilde x_k) dx^{k-1} \quad & \mbox{for }\, j \leq k \leq n-1 \\
0 \quad & \mbox{for }\, k < j
\end{cases}
\ee
\end{prop}
\begin{proof}
Recall $v = \varphi_n\circ \cdots \circ \varphi_1$ and
\be\label{eq:varphij}
\varphi_j(x_1,\cdots, x_n) = (x_1, \cdots, x_j
+ \zeta_j(x^{j-1})u(x_j,\widetilde x_j),\cdots, x_n).
\ee
We also note that we can write
$$
\int_{v(Q^n_{a;j})} g dy = \int_{Q^n_{a;j}} v^*(g\, dy)
$$
where
$$
v^*(g\, dy) = \varphi_1^*\circ \cdots \circ \varphi_n^*(g\, dy).
$$
Therefore to compute $d\Psi(\overline 0)\cdot X$, we need to first compute
the variation $\delta \varphi_j(X)$. But it is easy to see from
definition (\ref{eq:varphij}) of $\varphi_j$
\be\label{eq:deltaphij}
\delta \varphi_j(X) = (\zeta_j X_j) \frac{\del}{\del x_j}
\ee
and so
$$
(d\Psi(\overline 0)\cdot X)_j(a) = \sum_{k=1}^n \int_{Q^n_{a;j}}
\LL_{\delta \varphi_k(X)}(g\, dx) = \sum_{k=1}^n \int_{Q^n_{a;j}}
d(\delta \varphi_k(X)\rfloor (g\, dx)).
$$
On the other hand from the definition of $Q^n_{a;j}$, the
triangularity of $X$ and (\ref{eq:deltaphij}), the latter identity becomes
\beastar
(d\Psi(\overline 0)\cdot X)_j(a) & = & \sum_{k=j}^n \int_{Q^n_{a;j}}
d(\delta \varphi_k(X)\rfloor (g\, dx)) \\
& = & \sum_{k=j}^n \int_{Q^n_{a;j}} \frac{\del }{\del x_k}(g \zeta_k X_k)
dx_kdx^{n-1}_k \\
& = & \sum_{k=j}^n \int_{Q^{n-1}_{a;j,k}}g(x^{n-1}_{a;k}) \zeta_k(x^{k-1})
X_k(a_k,\widetilde x_k)
dx^{n-1}_k \\
& = & \sum_{k=j}^n \int_{Q^{n-k}_{a;j}} X_k(a_k,\widetilde x_k)
\Big(\int_{Q^{k-1}_{a;j}} \zeta_k(x^{k-1}) g(x^{n-1}_{a;k}) dx^{k-1}\Big)
d\widetilde x_k.
\eeastar
Here we define the $(n-1)$-vectors
\be
x^{n-1}_{a;k}  = (x^{k-1}, a_k, \widetilde x_k), \quad \mbox{for $k\geq j$}
\label{eq:yn-1aj}
\ee
and denote the volume element of any of $x^{n-1}_{a;k}$
by $dx^{n-1}_k$. Then the third equality above follows by
integration by parts over $x_k$.
This finishes the proof.
\end{proof}

Next we introduce the following function space which will be
essential for the later discussions :

\begin{defn}\label{Cvec1} We define
$$
C^{\vec 1}_{tri}(Q^n,\R^n)
$$
to be the set of continuous triangular maps
$f \in C^0_{tri}(Q^n,\R^n)$ whose components are given by the functions
$f_j : Q^{n-j} \to \R$ such that
$$
D^\alpha f : Q^{n-j} \to \R, \quad
Q^{n-j} = \{ (x_{j+1},\cdots, x_n) \mid
0 \leq x_l \leq 1, \, l = j+1, \cdots, n\}
$$
are continuous for any subset $\alpha \subset \{j+1, \cdots, n\}$.
Here $\vec 1$ stands for $\vec 1 = (1,\cdots, 1)$ and
$D^\alpha f$ for the partial derivative
with respect to the multi-index $\alpha$.
\end{defn}

It is easy to check that $C^{\vec 1}_{tri}(Q^n,\R^n)$ becomes
a Banach space if we equip it with a norm given by
$$
\|f\| = \max_{j =1, \cdots, n} \{\|f_j\|_{C^{\vec 1}}\}
$$
where $\|f_j\|_{C^{\vec 1}}$ is given by
\be\label{eq:|fj|Cvec1}
\|f_j\|_{C^{\vec 1}} = \max_{\alpha \subset \{j+1,\cdots, n\}}
|D^\alpha f_j|_{C^0}.
\ee
We recall that for any function $f \in C^{\vec 1}_{tri}(Q^n,\R^n)$
the partial derivatives $D^\alpha f$ does not depend on the
ordering of indices contained in the subset $\alpha \subset \{j+1,\cdots, n\}$
(See Theorem 7.3 \cite{lang}.)

With this preparation, we now prove
\begin{prop}\label{dPsi0vec1}
$d\Psi(\overline 0)$ continuously extends to a bounded linear operator
from $C^0_{tri}(Q^n,\R^n)$ to $C^{\vec 1}_{tri}(Q^n,\R^n)$ which is
bijective. Denote the extension again by
$$
d\Psi(\overline 0): C^0_{tri}(Q^n,\R^n) \to C^{\vec 1}_{tri}(Q^n,\R^n).
$$
In particular, it is invertible. We denote its inverse by
\be\label{eq:dPsiinverse}
(d\Psi(\overline 0))^{-1} : C^{\vec 1}_{tri}(Q^n,\R^n) \to
C^0_{tri}(Q^n,\R^n).
\ee
\end{prop}
\begin{proof}
From the matrix expression (\ref{eq:dPsi0}) of
$d\Psi(\overline 0) = \left(d\Psi(\overline 0)_{jk}\right)$,
we see that it becomes a triangular matrix and is represented by the integral pairing
with the functions $C_{jk}$ and manifestly
extends to an operator from $C^0(Q^{n-k},\R)$ to $C^{\vec 1}(Q^{n-j},\R)$.

And once we have proved the bijectivity of the bounded linear
operator $d\Psi(\overline 0)$, the open mapping theorem will imply that
the operator is invertible. Therefore it remains to prove bijectivity.

We start with the proof of injectivity. Suppose that $d\Psi(\overline 0)(X) = 0$
for $X = (X_1, \cdots, X_n) \in C^0_{tri}(Q^n,\R^n)$.
By (\ref{eq:dPsi0}) and (\ref{eq:Cjk}), $X$ satisfies
\be\label{eq:XC=0}
\sum_{k=j}^n \int_{Q^{n-k}_{a;j}} X_k(a_k,\widetilde x_k)
C_{jk}(a_j,\cdots, a_k,\widetilde x_k) d\widetilde x_k = 0
\ee
for all $j = 1, \cdots, n$. We will prove $X=0$ by a downward induction over $j$.
First consider the term for $j = n$.
In this case, this reduces to
$$
0= X_n(a_n) C_{nn}(a_n).
$$
Since $C_{nn}(a_n) =
\int_{Q^{n-1}_{a;n,n}}\zeta_n(x^{n-1})g(x^{n-1},a_n)dx^{n-1} > 0$,
we derive $X_n \equiv 0$.

Now suppose we have shown
$$
X_n = X_{n-1} = \cdots = X_{l+1} = 0
$$
and consider the equation $(d\Psi(\overline 0)(X))_l = 0$.
Under this assumption, (\ref{eq:XC=0}) for $j = l$ reduces to
$$
\int_{Q^{n-l}_{a;l}} X_l(a_l,\widetilde x_l) C_{ll}(a_l,\widetilde x_l) d\widetilde x_l = 0
$$
for all $a \in Q^n$. Differentiating this identity
with respect to $a_j$ successively for $j = l+1, \cdots, n$ at
the vector $a = (a_1, \cdots, a_n)$, we obtain
$$
0 = X_l(a_l,\widetilde a_l) C_{ll}(a_l,\widetilde a_l)
= X_l(a) C_{ll}(a_l,\widetilde a_l).
$$
Since $ C_{ll}(a_l,\widetilde a_l) > 0$, we obtain $X_l \equiv 0$
as before. This proves injectivity of $d\Psi(\overline 0)$.

Now we turn to surjectivity thereof. Let $Y = (Y_1, \cdots, Y_n)
\in C^{\vec 1}_{tri}(Q^n,\R^n)$ and consider the equation
$$
d\Psi(\overline 0)(X) = Y \quad \mbox{or equivalently }\,
(d\Psi(\overline 0)(X))_j = Y_j, \, j = 1, \cdots, n
$$
for $X \in C^0_{tri}(Q^n,\R^n)$. Again we solve this by downward induction
starting from $j = n$. For $j=n$, this reduces to
$$
X_n(a_n) C_{nn}(a_n) = Y_n(a_n)
$$
and so obtain $X_n(a_n) = Y_n(a_n) / C_{nn}(a_n)$. Now suppose that we have
solved for $j = n, \cdots, l+1$, and consider the equation
$(d\Psi(\overline 0)(X))_l = Y_l$. This equation becomes
\be\label{eq:XellYell}
\sum_{k=l}^n \int_{Q^{n-k}_{a;j}} X_k(a_k,\widetilde x_k)
C_{lk}(a_l,\cdots, a_k,\widetilde x_k) d\widetilde x_k = Y_l(a_l,\cdots,a_n).
\ee
Since $Y_l \in C^{\vec 1}(Q^{n-l},\R)$, we can differentiate
this equation with respect to $a_j$ successively over $j = l+1, \cdots, n$
and obtain
$$
X_l(a_l,\widetilde a_l) C_{ll}(a_l,\widetilde a_l)
+ \sum_{k = l+1}^n X_k(a_k,\widetilde a_k) \frac{\del C_{lk}}
{\del a_{l+1}\del a_{l+2}\cdots \del a_k} = \frac{\del Y_l}
{\del a_{l+1}\cdots\del a_n}
$$
by the triangularity of $X$ and $Y$. Since $C_{ll}(a_l,\widetilde a_l) > 0$,
we obtain
\beastar
X_l(a) & = & X_l(a_l,\widetilde a_l)\\
& = & \frac{1}{C_{ll}(a_l,\widetilde a_l)}
\left(\frac{\del Y_l}{\del a_{l+1}\cdots\del a_n}(a)
- \sum_{k = l+1}^n X_k(a) \frac{\del C_{lk}}
{\del a_{l+1}\del a_{l+2}\cdots \del a_k}(a)\right).
\eeastar
We note that by the induction hypothesis, the right hand side
is already determined.
Since $Y_l \in C^{\vec 1}(Q^{n-l},\R)$ and $C_{lk}$ are smooth,
the right hand side is continuous
and hence lies in $C^0_{tri}(Q^n,\R^n)$. This finishes the induction step
and so solves the equation
$d\Psi(\overline 0)(X) = Y$ for any $Y \in C^{\vec 1}_{tri}(Q^n,\R^n)$ and so
finishes the proof of surjectivity. Hence the proof.
\end{proof}

The following proposition is a crucial ingredient which
saves us from doing derivative estimates for the nonlinear
terms in section \ref{sec:higher}.

\begin{prop}\label{C0extend}
The operator $(d\Psi(\overline 0))^{-1}$ given in (\ref{eq:dPsiinverse})
continuously extends to a bounded linear operator
$$
K_g : C^0_{tri}(Q^n,\R^n) \to C^0_{tri}(Q^n,\R^n).
$$
\end{prop}
\begin{proof} We go back to the surjectivity proof of Proposition
\ref{dPsi0vec1}. It will be enough to prove that there
exists a constant $M_g > 0$ such that the unique solution $X$ for
$d\Psi(\overline 0)(X) = Y$ satisfies
\be\label{eq:XC0}
|X|_{C^0} \leq M_g |Y|_{C^0}
\ee
for any given $Y \in C^{\vec 1}_{tri}(Q^n,\R^n)$. We prove this again
by the downward induction.

For $j = n$, we have $X_n(a_n) C_{nn}(a_n) = Y_n(a_n)$ and hence
$$
|X_n(a_n)| \leq \frac{|Y_n(a_n)|}{\min_{a_n}C_{nn}(a_n)}
$$
recalling $C_{nn}(a_n) \neq 0$ from (\ref{eq:Cjk}). But we have
$$
C_{nn}(a_n) = \int_{Q^{n-1}_{a;n}}\zeta_n(x^{n-1})g(x^{n-1},a_n)dx^{n-1}
$$
from (\ref{eq:Cjk}). In particular, we have
$$
\min_{a_n}C_{nn}(a_n) \geq \min g \int_{Q^{n-1}_{a;n}}\zeta_n(x^{n-1})dx^{n-1}
= \min g
$$
where we use (\ref{eq:zetas}) for the equality.
Therefore we have proved
\be\label{eq:Xnbound}
|X_n|_{C^0} \leq \frac{1}{\min g} |Y_n|_{C^0}
\ee
for all $l+1 \leq j \leq n$. Now as the induction hypothesis,
suppose that there exists a constant $M_{\ell+1}$ such that
\be\label{eq:Xell+1bound}
|X_j|_{C^0} \leq M_{\ell +1} \max_{ \ell +1 \leq j \leq n} |Y_j|_{C^0}
\ee
for all $l+1 \leq j \leq n$.
We rewrite (\ref{eq:XellYell}) into
\bea\label{eq:Xell=Yell}
\int_{Q^{n-l}_{a;l}} X_l(a_l,\widetilde x_l)
C_{ll}(a_l,\widetilde x_l) d\widetilde x_l =
Y_l(a_l,\cdots,a_n)\nonumber \\
\quad - \sum_{k=l+1}^n \int_{Q^{n-k}_{a;j}} X_k(a_k,\widetilde x_k)
C_{lk}(a_l,\cdots, a_k,\widetilde x_k) d\widetilde x_k.
\eea
By the induction hypothesis the sum in the right hand side can be
estimated as
\beastar
&{}& \left|\sum_{k=l+1}^n \int_{Q^{n-k}_{a;j}} X_k(a_k,\widetilde x_k)
C_{lk}(a_l,\cdots, a_k,\widetilde x_k) d\widetilde x_k\right| \\
& \leq & \int_{Q^{n-k}_{a;j}}M_{\ell +1}\max_{ \ell +1 \leq j \leq n}
\{|Y_j|_{C^0}\}C_{lk}(a_l,\cdots, a_k,\widetilde x_k) d\widetilde x_k\\
& \leq & M_{\ell +1}\max_{ \ell +1 \leq j \leq n}
\{|Y_j|_{C^0}\}\int_{Q^{n-k}_{a;j}}C_{lk}(a_l,\cdots, a_k,\widetilde x_k)
d\widetilde x_k.
\eeastar
But we derive
\beastar
\int_{Q^{n-k}_{a;j}}C_{lk}(a_l,\cdots, a_k,\widetilde x_k)
d\widetilde x_k & = & \int_{Q^{n-k}_{a;\ell}}
\left(\int_{Q^{k-1}_{a;\ell}}\zeta_k(x^{k-1})g(x^{k-1},a_k,\widetilde x_k)dx^{k-1}
\right) d\widetilde x_k \\
&\leq & \max g \left(\int_{Q^{k-1}_{a;\ell}}\zeta_k(x^{k-1})dx^{k-1}
\right) \leq \max g
\eeastar
again using (\ref{eq:zetas}).
Hence we have obtained
\beastar
&{}&\left|\sum_{k=l+1}^n \int_{Q^{n-k}_{a;j}} X_k(a_k,\widetilde x_k)
C_{lk}(a_l,\cdots, a_k,\widetilde x_k) d\widetilde x_k \right| \\
&\leq &(n -\ell -1) M_{\ell +1}\max g \max_{\ell +1 \leq j \leq n}
\{|Y_j|_{C^0}\}.
\eeastar
Substituting this into (\ref{eq:Xell=Yell}), we obtain
\beastar
&{}& \left|\int_{Q^{n-l}_{a;l}} X_l(a_l,\widetilde x_l)
C_{ll}(a_l,\widetilde x_l) d\widetilde x_l\right|\\
& \leq &
|Y_l|_{C^0} + (n-\ell -1) M_{\ell +1}\max g \max_{\ell +1 \leq j \leq n}
\{|Y_j|_{C^0}\} \\
& \leq & (n-\ell) \max\{M_{\ell},1\} \max g \max_{\ell \leq j
\leq n} \{|Y_j|_{C^0}\}\}.
\eeastar

On the other hand, the left hand
side can be estimated from below
\beastar
\left|\int_{Q^{n-l}_{a;l}} X_l(a_l,\widetilde x_l)
C_{lk}(a_l,\widetilde x_l) d\widetilde x_l\right|
&\geq &|X_l(a_l,\widetilde x_l)| \min C_{ll} \\
& \geq & |X_l(a_l,\widetilde x_l)| \min g.
\eeastar
Combining the last two inequalities, we have obtained
$$
|X_l(a_l,\widetilde x_l)| \leq (n-\ell)\frac{\max\{M_{\ell+1},1\}
\max g \max_{\ell \leq j \leq n} \{|Y_j|_{C^0}\}}{\min g}.
$$
By defining
$$
M_\ell = (n-\ell) \frac{\max\{M_{\ell+1},1\} \max g}{\min g}
$$
we have finished the induction step and hence the proof of (\ref{eq:XC0}).

In fact the above proof shows that $M_g$ can be chosen to be
\be\label{eq:Mg} M_g
= n!\max\left\{\frac{1}{\min g}
\left(\frac{\max g}{\min g}\right)^{n-1}, 1\right\} \ee and hence we
have $\|K_g\| \leq M_g$. This finishes the proof.
\end{proof}

We recall that the constant $M_g$ does not increase under the
reduction process to a smaller cubes by the reasons mentioned in
Remark \ref{Qn-X}.


\section{$C^0$-coerciveness of Darcorogna-Moser's solutions}
\label{sec:functional}

We denote the operator norm of the bounded linear operator
$K_g$ given in Proposition \ref{C0extend} by $\|K_g\|$ which
has the bound
\be\label{eq:Modg}
\|K_g\| \leq M_g
\ee
where $M_g$ is the constant given in (\ref{eq:Mg}).
From the explicit formula of $M_g$, it follows that
$M_g$ depends only on $g$ and is continuous on $g$
in $C^0$-topology.

In this section, all the norms $|\cdot|$ below will
denote the $C^0$-norms.

We write (\ref{eq:-dPsi0}) in the following form
\be\label{eq:uXiu}
u = \Xi(u)
\ee
where $\Xi$ is the map from $B \to B$ defined by
\be\label{eq:Xiu}
\Xi(u) = - (d\Psi(\overline 0))^{-1}(\Psi(\overline 0) + N(u)).
\ee
Here we would like to note from (\ref{eq:Psia0})
that $\Psi(\overline 0)$ lies in $C^{\vec 1}(Q^n,\R^n)$.
On the other hand, we can rewrite $N(u)$
\bea\label{eq:N(u)}
N_j(u) & = &\Psi_j(u) - \Psi_j(\overline 0) - (d\Psi(\overline 0)\cdot u)_j
\nonumber\\
& = & \int_{R^n_{au;j}} g(y) dy -
\int_{Q^n_{a;j}} g(y) dy \nonumber\\
&{}& \hskip0.2in -\sum_{k=j}^n \int_{Q^{n-k}_{a;j}} u_k(a_k,\widetilde x_k)
\Big(\int_{Q^{k-1}_{a;j,k}} \zeta_j(x^{k-1}) g(x^{n-1}_{a;k}) dx^{k-1}\Big)
d\widetilde x_k.
\eea
From this, it follows that $N(u)$ also lies in $C^{\vec 1}_{tri}(Q^n,\R^n)$
if $u$ is smooth as for $u = u_{DM}$. Therefore $\Psi(\overline 0) + N(u)$
lies in the domain of $(d\Psi(\overline 0))^{-1}$ and hence
the expression (\ref{eq:Xiu}) is well-defined for $u = u_{DM}$.

We derive from (\ref{eq:uXiu}), (\ref{eq:Xiu})
\be\label{eq:|u|<|Xiu|}
|u| \leq |(d\Psi(\overline 0))^{-1}(\Psi(\overline 0) + N(u))|
\leq M_g (|\Psi(\overline 0)| + |N(u)|).
\ee
We now estimate $|\Psi(\overline 0)|$ and  $|N(u)|$ separately.

We start with the following

\begin{prop}\label{|Psi0|} We have
\be
|\Psi(\overline 0)| \leq d_\CM(m_f,m_g)
\ee
where $m_f := \mu_{(fdx)}$ is the measure associated to
the volume form $fdx$ and similarly for $m_g$.
\end{prop}
\begin{proof}
From (\ref{eq:Psia0}), we have for the $j$-th component of
the vector $\Psi(\overline 0)(a)$
\be\label{eq:Psi0a}
\Psi(\overline 0)_j(a) = G_j(a,0) - F_j(a) =
\int_{Q^n_{a;j}} g(y) dy - \int_{Q^n_{a;j}} f(x) dx.
\ee
On the other hand $R_{a;j} = Q_{a;j}$ for $u = \overline 0$. Therefore
from the definition (\ref{eq:Lidb}) of the metric $d_\MM=\operatorname{Lid}_1$, we
have derived the upper-bound for the `zero-order term'
\be\label{eq:|Psi0|}
|\Psi(\overline 0)| \leq d_\MM(m_f,m_g)
\ee
where we use the fact that the integral (\ref{eq:Psi0a}) corresponds to
$$
\int \chi_{Q^n_{a;j}}\,dm_g
- \int \chi_{Q^n_{a;j}}\, dm_f
$$
which is obtained by taking the characteristic function $\chi_{Q^n_{a;j}}$
of $Q^n_{a;j}$ as the test function in (\ref{eq:Lidb}) for $b=1$.
\end{proof}

Proposition \ref{|Psi0|} is a place where \emph{the weak topology of
measures enters in our proof of $C^0$-coerciveness} (\ref{eq:dMMX})
in Theorem II. The other such places appearing later will be similar
to this one.

Next we do estimates of $N(u)$. For this purpose, we introduce the
constant \be\label{eq:e1zeta} \e_1(\zeta) = \max_{1\leq k \leq
n}\left( \int_{Q^{k-1}}|1 - \zeta_k(y^{k-1})| dy^{k-1}\right). \ee

See the end of section \ref{sec:higher} for our motivation for
considering this constant where it appears in middle of the main
technical estimates. We like to emphasize that \emph{this constant
can be made as small as we want by approximating $\zeta$
$L^1$-close to the function $1$, once $g, \, f$ are given}.
In particular, we may assume
\be\label{eq:e1<}
\e_1(\zeta) < d_\MM(m_f,m_g).
\ee

Next using the continuity of $g$ and compactness of $Q^n$, we
have the Lipschiz bound
\be\label{eq:Lipschiz}
|g(x) - g(y)| \leq L_g \cdot |x-y|
\ee
for a constant $L_g > 0$ depending only on $g$. In fact, $L_g$ is
nothing but the \emph{modulus of continuity} of $g$.

The following is a key lemma whose proof we postpone until
section \ref{sec:higher} because the proof is rather long and
complicated. The main reason behind the presence of this
kind of estimates is that $N(u)$ is
the higher order term in the Taylor expansion of $\Psi$.
However, since we need to
know the precise form of the inequality with respect to $g$, we need
to carry out rather delicate estimates.

\begin{lem}\label{lemma-N(u)} Define
$$
C_4 = C_4(g) =  \max\{8\max g, 4\}
$$
and let $u = u_{DM}$ be a DM-solution.
Then we have the inequality
\be\label{eq:NuC3}
|N(u)| \leq C_4 \cdot\left(|\zeta\cdot u|+ \e_1(\zeta)
+ L_g |\zeta\cdot u|\right)|u|,
\ee
where we denote $\zeta\cdot u:= (u_1,\zeta_2u_2, \cdots, \zeta_nu_n)$
\end{lem}

Combining Proposition (\ref{eq:|Psi0|}), (\ref{eq:NuC3}),
and (\ref{eq:|u|<|Xiu|}), we obtain
\be\label{eq:|u|leqMg}
|u| \leq M_g (d_\MM(m_f,m_g) +
C_4 \cdot\left(|\zeta\cdot u|+ \e_1(\zeta)
+ L_g |\zeta\cdot u|\right)|u|).
\ee

We can choose
the functions $\zeta=\{\zeta_k\}_{k =2, \cdots, n}$ so that
$$
\e_1(\zeta) < d_\MM(m_f,m_g)
$$
as mentioned in (\ref{eq:e1<}). Substituting this into and
rewriting (\ref{eq:|u|leqMg}), we obtain
\be\label{eq:1-C4}
\Big(1 - C_4 \cdot(|\zeta \cdot u| +d_\MM(m_f,m_g)+ L_g
\cdot |\zeta \cdot u|)\Big)|u| \leq M_g \cdot d_\MM(m_f,m_g).
\ee

At this stage, we recall that a DM-solution has the form
$$
v = (v_1,v_2, \cdots, v_n)
$$
with $v_j(x) = x_j + \zeta_j(v^{j-1}(x)) u_j(x_j, \widetilde x_j)$
and maps $Q^{n}$ into $Q^{n}$. In particular, we have
$$
|\zeta_ju_j| = \max_{x \in Q^n}|\zeta_j(x) u_j(x)| \leq 2.
$$

We would also like to emphasize that the constants $C_4(g)$
and $L_g$ depend only on $g$ but not on $f$, except in the loose way
mentioned in (\ref{eq:e0}). Therefore Dacorogna-Moser's
construction of solution can be
equally carried out for the maps defined on the cube $Q^n(K)$ with
any length $0 < K \leq 1$ of its sides with the same constants
$C_4$ and $L_g$. In that case, all DM-solutions
on $Q^n(K)$ will satisfy
$$
|\zeta \cdot u|\leq 2 K
$$
because $u$ maps $Q^n(K)$ to $Q^n(K)$ in that case.

Therefore if we set $K = 2^{-N_0}$ and fix $N_0= N_0(g) \in \N$ such that
\be\label{eq:C4}
C_4 \cdot (1 + L_g)
\cdot 2^{-(N_0-1)}< \frac{1}{4}
\ee
and consider a DM-solution on the cube $Q^n(K)$,
we will have
\be\label{eq:C4cdot}
C_4 \cdot (2^{-(N_0-1)} + d_\MM(m_f,m_g)+ L_g
\cdot 2^{-(N_0-1)}) < \frac{1}{2}
\ee
for any $f$ such that
$$
C_4 d_\MM(m_f,m_g) < \frac{1}{4}.
$$
We recall from Remark \ref{Qn-X} that this inequality will be
achieved by considering $\mu_\sigma, \, \mu_\tau$ with
$d_{\MM}(\mu_\sigma,\mu_\tau) \to 0$ on the original space $X$
given in Theorem II.

Then (\ref{eq:1-C4}) and (\ref{eq:C4cdot}) imply
$$
|u| \leq 2 M_g\cdot d_\MM(m_f,m_g).
$$
Here we recall from (\ref{eq:m0M0}) that $M_g$
depends only on the originally given
function $g$ defined on the unit cube $Q^n$.

Now by decomposing $Q^n$ into cubes of size $2^{-N_0}$,
and applying this inequality uniformly over to each of
the cubes, we obtain
the following proposition. Here $N_0$ is the integer
chosen as in (\ref{eq:C4}), which depends only on the
originally given function $g$ defined on $Q^n$.

\begin{prop}\label{fixedXi}
Let $g$ be a positive continuous function and
denote by $m_g = \mu_{(g dx)}$ the associated measure
on $Q^n$. Consider the Darcorogna-Moser's
solution $u = u_{DM}$ corresponding to $f$ satisfying
the hypotheses in Theorem 3.2.
Then there exists a continuous function $r = r(t;g)$
of $t$, depending only on $g$, such that $r \to 0$ as $t \to 0$
for which the following holds :
\par
\be
|u|_{C^0} \leq r(d_\MM(m_f,m_g);g).
\ee
\end{prop}

To wrap-up the proof of statement (1) of Theorem \ref{c0est}, we need
to estimate $d_{C^0}(\psi^{-1}_2,id)$ for $\psi_2 = v$.

For the estimate of  $d_{C^0}(\psi^{-1}_2,id)$,
we derive $d_{C^0}(\psi^{-1}_2,id) = d_{C^0}(id, \psi_2)$. For we have
$$
d_{C^0}(\psi^{-1}_2,id) = \max_{x \in Q^n}d(\psi^{-1}_2(x),x)
= \max_{x \in Q^n}d(\psi^{-1}_2(x),\psi_2(\psi^{-1}_2(x)) \leq d_{C_0}(id,
\psi_2)
$$
and prove the opposite inequality in the same way. Therefore we obtain
$$
\bar d(\psi_2,id) \leq d_{C^0}(\psi_2,id) + d_{C^0}(\psi_2^{-1},id)
\leq 2 r(d_\MM(m_f,m_g);g).
$$
This finishes the proof of (1) of Theorem \ref{c0est} and
hence the proof of Theorem \ref{c0est}.

\section{Proof of Theorem II}
\label{sec:theoremII}

With Proposition \ref{local} and Theorem \ref{c0est} in our hand,
we now give the proof of Theorem II. We will imitate
Moser's argument \cite{moser} but with some additional arguments
needed to establish the $C^0$-coerciveness.

Let $\tau = \Omega$ be a volume form on $X$, $f > 0$ be a positive
function on $X$ and $\sigma = f\Omega$. We choose an open covering
$\UU = \{U_j\}$ of $X$ and denote by $m=m_\UU$ the cardinality of $\UU$.

Consider the functions $f(t;\cdot)$ defined by
$$
f(t;p) = 1 + \sum_{j=0}^m t_j g_j(p), \quad t = (t_0, \cdots, t_m).
$$
For $t = (0,\cdots, 0)$, one has $f(t;p) \equiv 0$ and for
$t = (1,\cdots,1)$, we have $f(t;\cdot) = f$. By construction,
we also have $\int f(t;\cdot) \Omega = \int \Omega$ for all $t$.
We can connect two corners $(0,\cdots, 0)$ and $(1,\cdots,1)$
of the cube by going along $m+1$ edges. If $t',\, t''$ represent
the endpoints of such an edge, one sees that
$$
f(t'';\cdot) - f(t';\cdot)
$$
has support in one patch, say $U_{\{t',t''\}}$. Without loss of any generality,
we may parameterize
$$
U_{\{t',t''\}} \cong (-\eta,1+\eta)^n
$$
for some $\eta > 0$ and
$$
\text{supp}(f(t'';\cdot) - f(t';\cdot)) \subset Q^n(1-\eta)
$$
under the parametrization. If we write
$$
\Omega_t = f(t;\cdot) \Omega
$$
one sees that $\Omega_{t''} = h \Omega_{t'}$ where
$h = f(t'';\cdot)/f(t';\cdot)$ is different from 1 in $U_{\{t',t''\}}$
only and $h \equiv 1$ on $Q^n(1+\eta) \setminus Q^n(1-\eta)$.

Once we have made the choice of such a covering $U_j$,
we consider the family
$$
(t,s) \in [0,1]^{n+1} \times [0,1] \mapsto f(t,s)
$$
where $f(t,s)$ is the function on $X$ defined by
\be\label{eq:ftsp}
f(t,s;p) = 1 + s \sum_{j=0}^m t_j g_j(p).
\ee
We partition $s \in [0,1]$ into a partition
$$
P: s_0 = 0 < s_1 < \cdots < s_N.
$$
By choosing $P$ with $\text{mesh}(P)$ sufficiently small, we
can make
\be\label{eq:meshP}
d_\MM(m_{f(t',s_i)},m_{f(t'',s_i)})
\ee
as small as we want uniformly over $(t',t'')$. Therefore we
will assume that $d_\MM(\mu_{\Omega},\mu_{f\Omega})$ is so small that we can
apply Proposition \ref{fixedXi}.

We order the set $U_j$ so that $U_j = U_{\{t_j,t_{j+1}\}}$
where $U_{\{t_j,t_{j+1}\}}$ is the patch corresponding to
$(t',t'') = (t_j,t_{j+1})$.
Now applying Theorem \ref{c0est} to each patch $U_j$,
we have constructed a sequence of diffeomorphisms $\phi_j: X \to X$
such that
\begin{enumerate}
\item $\phi_{j+1} \circ \phi_j^{-1}$ has support in
$Q^n(1-\eta) \subset U_j\cong Q^n(1+\eta)$
\item $\Omega_{t_{j+1}} = \phi_j^*\Omega_{t_j}$ or equivalently
$(\phi_i)_*\mu_{j+1} = \mu_j$
where $\mu_j = \mu_{\Omega_{t_j}}$.
Here we denote by $t_j$ the $j$-th vertex in the above chosen
edge path from $(0,\cdots,0)$ to $(1,\cdots,1)$.
\end{enumerate}
Then the diffeomorphism $\psi_2 = \phi_m: X \to X$ satisfies
$\psi_2^*\Omega = f\Omega$.

It remains to estimate $\bar d(\psi_2,id)$.
Since $d_{C^0}(\psi_2^{-1},id) = d_{C^0}(id, \psi_2)$, it is enough
to estimate $d_{C^0}(\psi_2,id)$.
Denote the above coordinate patch map
$$
\psi_j : U_j \to (-\eta,1+\eta)^n.
$$
Then the above diffeomorphism $\phi_{j+1}\circ \phi_j^{-1}$
is given by the conjugation
$$
\phi_{j+1}\circ \phi_j^{-1} = \psi_j^{-1} \circ \varphi_j \circ
\psi_j
$$
where $\varphi_j: Q^n(1+\eta) \to Q^n(1+\eta)$ is the diffeomorphism
constructed in section \ref{sec:functional} corresponding to
the forms
$$
(f_{j+1}\circ \psi_{j+1}^{-1})\, dx, \quad (f_j \circ \psi_j^{-1})\, dx.
$$
Because we will use the $C^0$-norm in different spaces, we will specify
the space where the $C^0$-norm is taken below when we need to specify
the space. We have from Proposition \ref{fixedXi}
\be\label{eq:C0Qj}
d_{(C^0, Q^n)}(\varphi_j,id) \leq
C_5(d_\MM(m_{j+1},m_j);f_j \circ \psi_j^{-1})
\ee
where $m_j$ is the measure associated to the form
$(f_j \circ \psi_j^{-1})\, dx$. We recall that the finite family of functions
$f_j\circ \psi_j^{-1}$ are determined by the original function $f$, the covering
$\UU$ and the coordinate charts $\psi_j$.
Since we do not change but fix them in the course of proof,
we may ignore this dependence of $C_5$ on $f_j\circ\psi_j^{-1}$.

We note that
$$
d_{(C^0, X)}(\phi_{j+1}\circ \phi_j,id) = d_{C^0}
(\psi_j^{-1}\circ \varphi_j \circ \psi_j,id) \to 0
$$
as
$d_{(C^0,Q^n(1+\eta))}(\varphi_j,id) \to 0$ and
\be\label{eq:C0Xpsi2id}
d_{(C^0, X)}(\psi_2,id) \leq \sum_{j=0}^{m-1}
d_{(C^0, X)}(\phi_{j+1},\phi_{j}).
\ee
On the other hand, we have
\bea
d_{(C^0, X)}(\phi_{j+1},\phi_{j}) & = & d_{(C^0, X)}((\phi_{j+1}\circ
\phi_j^{-1})\circ \phi_j, \phi_j) \nonumber \\
& = & d_{(C^0, X)}(\phi_{j+1}\circ \phi_j^{-1}), id)\nonumber \\
& = & d_{(C^0, X)}(\psi_j^{-1}\circ \varphi_j \circ \psi_j,id).
\label{eq:C0Xj}
\eea
Since the integer $m = m_\UU < \infty$ depends only on
$\UU$ but not on $f$ or $j$'s, it follows that as $d_\MM(\mu_{f\Omega},\mu_\Omega) \to 0$,
$\bar d(\psi_2,id) \to 0$. This finishes the proof of Theorem II
except its parameterized version.

For the parameterized version, we recall that $\tau = \Omega$ is fixed
and that $\sigma_s = f_s \Omega$ for a given a smooth family of functions $f_s$
for $s \in [0,1]$, for which $s \mapsto \mu_{(f_s\Omega)}$
is continuous in $\MM(X)$. Note that the above mentioned covering $\UU = \{U_j\}$
in section \ref{sec:cube} does not depend on the functions $f_s$
and so can be fixed for all $s \in [0,1]$.
This and the compactness of $[0,1]$ enable us to reduce the problem
to the parameterized version of Theorem \ref{c0est}
on the cube for a fixed $g$ but varying $f_s$ in a way that
$s \mapsto \mu_{(f_sdx)}=m_{f_s}$ is continuous in $\MM(Q^n)$.
Since all the constants appearing in section \ref{sec:functional}
depend continuously on $d_\MM(m_f,m_g)$, we can uniformly apply
Dacorogna-Moser's construction  to produce an isotopy $s \mapsto
\psi_{2,s}$ of diffeomorphisms that is continuous in compact
open topology of $Diff(X)$ and satisfies $(\psi_{2,s}^{-1})^*\Omega
= f_s\Omega$ which is equivalent to $\Omega
= \psi_{2,s}^*(f_s\Omega)$.

One particular remark on the choice of
the constant $\e_0 = \e_0(g,\zeta)$ in our construction on the cube
$Q^n$ is in order for the parameterized case. For the given isotopy
$\FF=\{f_s\}_{0 \leq s \leq 1}$ on $M$, we can reduce the problem
to the cube so that
$$
\operatorname{supp}(f_s - g ) \subset Q^n(1-\eta)
$$
for all $s \in [0,1]$. (See section \ref{sec:cube}.)
Then we choose $\e_0^{\FF}=\e_0(\FF,\zeta)$ by
$$
\e_0^{\FF} = \min_{s\in [0,1]}\e_0(f_s,\zeta)
$$
which can be made close to 0 uniformly over $s$
by choosing the family $\zeta = \{\zeta_s\}_{s=2}^n$
of cut-off functions $\zeta$ as in (\ref{eq:zetas})
suitably.

To improve the regularity of the parameterized solutions,
we use the a priori $C^k$ estimate provided in Proposition \ref{Ckestimate}.
This finishes the proof of Theorem II

\section{Proof of Theorem $\mbox{\rm I}'$}
\label{sec:TheoremI}

In this section, we finish the proof of
Theorem $\mbox{\rm I}'$ following the scheme outlined in the introduction.

Let $h \in M[\Sigma,\Omega]$ and $\e > 0$ be given.
By the smoothing theorem
(see Theorem 6.3 \cite{munkres:smoothing} for example), we can choose a diffeomorphism
$\psi_1$ such that
\be\label{eq:dC0}
\bar d(h,\psi_1) \leq \frac{\e} 2.
\ee
This diffeomorphism is {\it not} necessarily
area preserving. We therefore modify $\psi_1$ into an area
preserving diffeomorphism by composing it with another diffeomorphism
$\psi_2:\Sigma \to \Sigma$ that is $C^0$-close to the identity.

It follows from Proposition \ref{inMMX} that (\ref{eq:dC0}) also implies that
the measures associated to $\mu_{\Omega}=h_*(\mu_\Omega)$ and
$(\psi_1^{-1})_*(\mu_\Omega)= \mu_{((\psi_1)^*\Omega)}$
can be made arbitrarily close in the weak topology
of measures.

We note that $\int_\Sigma (\psi_1)^*\Omega
= \int_\Sigma \Omega$. Therefore
applying Theorem II to the forms
$$
\sigma = \Omega, \quad \tau = (\psi_1^{-1})^*\Omega, \quad
\mbox{with }\, \lambda = 1
$$
we obtain a diffeomorphism $\psi_2$ such that
$$
\psi_2^*\Omega = (\psi_1^{-1})^*\Omega
$$
and
\be\label{eq:dC0psi2id}
\bar d(\psi_2,id) \to 0 \quad \mbox{as $d_\MM(\mu_{((\psi_1)^*\Omega)},
\mu_{\Omega}) \to 0$}.
\ee

We set $\phi = \psi_2 \circ \psi_1$.
Then $\phi$ is an $\mu_\Omega$-area preserving diffeomorphism and we have
\be\label{eq:c0psih}
\bar d(\phi,h) \leq \bar d(\psi_2 \circ \psi_1,\psi_1)
+ \bar d(\psi_1,h)
\ee
by the triangle inequality. But since $X$ is compact and $\psi_1$ is a
diffeomorphism, we can make $\bar d(\psi_2\circ \psi_1,\psi_1)$ as small
as we want by having $d(\psi_2,id)$ sufficiently small.
However (\ref{eq:dC0psi2id}) implies that $d(\psi_2,id)$
can be made as small as we want if we can let $d_{\MM}(\mu_{(\psi_1^*\Omega)}, \mu_\Omega)
= d_{\MM}((\psi_1^{-1})_*\mu_{\Omega}, \mu_\Omega)$ arbitrarily small.
And the latter can be achieved by Proposition
\ref{inMMX} if we choose the initial diffeomorphism
$\psi_1$ sufficiently $C^0$ close to $h$.
Combining these with (\ref{eq:c0psih}), we
can make
$$
\bar d(\phi,h) \leq \bar d(\psi_2 \circ \psi_1,\psi_1) + d(\psi_1,h)
\leq \frac{\e}{2} + \frac{\e}{2} = \e
$$
if we choose the initial smooth approximation $\psi_1$
sufficiently $C^0$ close to $h$.
This finishes the proof of Theorem $\mbox{\rm I}'$ (1).

Finally when we are given an isotopy $\{h_t\}_{0 \leq t\leq 1}$
of homeomorphisms $h_t: \Sigma \to \Sigma$, we apply  the isotopy version
of smoothing theorem (see Theorem 6.3 \cite{munkres:smoothing})
to obtain an isotopy of diffeomorphisms $\{\psi_{1,t}\}_{0 \leq t\leq 1}$
so that $\bar d(\psi_{1,t},h_t)$ can be made as small as we want
uniformly over $0 \leq t \leq 1$.

Then the isotopy of forms
$$
t \in [0,1] \mapsto (\psi_{1,t})^*\Omega
$$
is continuous in the sense mentioned in Theorem II. Therefore
we can apply the parameterized version of Theorem II to produce another
isotopy $\{\psi_{2,t}\}_{0 \leq t\leq 1}$ so that
\begin{enumerate}
\item $\psi_{1,t}^*(\psi_{2,t}^*\Omega) = \Omega$ for all $0 \leq t\leq 1$.
\item The isotopy $t \mapsto \psi_{2,t}$ is continuous in compact open topology.
\item $\bar d(\psi_{2,t},id)$ is as small as we want uniformly over $t \in [0,1]$.
\end{enumerate}
Now the composed isotopy defined by $\psi_t = \psi_{2,t}\circ \psi_{1,t}$
will do our purpose. This finishes the proof of Theorem $\mbox{\rm I}'$.

\begin{rem}
Here we would like to point out that we should apply our construction
to the \emph{fixed} form $\Omega$, not to the \emph{varying} form $\psi_1^*\Omega$.
In this way, the dependence on $g$ appearing in all the constants
in section \ref{sec:functional} become irrelevant in our construction
because this function $g$ will be fixed throughout the construction.
This is the reason why we consider the pair of forms
$$
\sigma = \Omega, \quad \tau = (\psi_1^{-1})^*\Omega
$$
instead of the more naturally looking choice of
$$
\sigma = (\psi_1)^* \Omega, \quad \tau = \Omega.
$$
\end{rem}

\section{Estimates of the higher order terms}
\label{sec:higher}

In this section, we prove Lemma \ref{lemma-N(u)}. We would like to
note that the nonlinear terms depends only on $g$. Dependence
on $f$ occurs only through the constants $\e_0(\zeta)$ and
$\e_1(\zeta)$ which involves the choice of cut-off functions $\zeta$.
However these constants can be made as small as we want, e.g.,
smaller than $d_\MM(m_f,m_g)$ as in (\ref{eq:e0}) and (\ref{eq:e1zeta})
by choosing $\zeta$ appropriately, once $f$ is given.

We start with the statement (1) of Lemma 5.2. We recall from (\ref{eq:N(u)})
\beastar
N_j(u)
& = & \int_{R^n_{au;j}} g(y) dy -
\int_{Q^n_{a;j}} g(y) dy \nonumber\\
&{}& \hskip0.2in -\sum_{k=j}^n \int_{Q^{n-k}_{a;j}} u_k(a_k,\widetilde x_k)
\Big(\int_{Q^{k-1}_{a;j,k}} \zeta_j(x^{k-1}) g(x^{n-1}_{a;k}) dx^{k-1}\Big)
d\widetilde x_k.
\eeastar
Now to estimate (\ref{eq:N(u)}), we rewrite $v(Q^n_{a;j})$ as
$$
v(Q^n_{a;j}) = \varphi_n\circ\varphi_{n-1} \circ\cdots
\circ \varphi_j(Q^n_{a;j})
$$
by the factorization of
$v = \varphi_n\circ\varphi_{n-1} \circ\cdots \circ \varphi_1$.
Here we also use the identity
$$
\varphi_{j-1}\circ \cdots \circ \varphi_1(Q^n_{a;j})
=Q^n_{a;j}
$$
which follows from the definitions of $\varphi_k$ and $Q^n_{a;j}$.
Motivated by this,
for each $j-1 \leq l\leq n+1$, we define
$$
R^n_{au,l;j} =  \varphi_n \circ \cdots \circ \varphi_{l}(Q^n_{a;j}).
$$
We note that $R^n_{au,j;j} = R^n_{au;j}$
and set $R^n_{au,n+1;j} = Q^n_{a;j}$.

From now on, we will switch the variable $y$ with $x$
and use $x$ instead of $y$ for the rest of the proof.

With these definitions, we can telescope and rewrite
$$
\int_{R^n_{au;j}} g(x) dx - \int_{Q^n_{a;j}} g(x) dx
= \sum_{l=j}^{n} \left(\int_{R^n_{au,l;j}} g(x) dx
- \int_{R^n_{au,l+1;j}} g(x) dx\right).
$$
On the other hand, we can easily check
\bea\label{eq:uka}
\int_{Q^{n-k}_{a;j}} u_k(a_k,\widetilde x_k)
\Big(\int_{Q^{k-1}_{a;j,k}} \zeta_j(x^{k-1}) g(x^{n-1}_{a;k}) dx^{k-1}\Big)
d\widetilde x_k = \nonumber\\
\int_{Q^n_{a;j,k}\setminus Q^n_{a;j}} \zeta_k(x^{k-1}) g(x^{n-1}_{a;k})dx
-\int_{Q^n_{a;j}\setminus Q^n_{a;j,k}} \zeta_k(x^{k-1})g(x^{n-1}_{a;k})dx
\eea
where $Q^n_{a;j,k}$ is defined by
\bea\label{eq:Qnajk}
Q^n_{a;j,k} =  \{ x \in Q^n & \mid & 0 \leq x_i \leq 1 \mbox{ for } i = 1, \cdots, j-1,
\nonumber\\
& {} &  0 \leq x_k \leq a_k + u_k(a_k,\widetilde x_k) \nonumber\\
& {} &  0 \leq x_i \leq a_i \mbox{ for } i \geq j \, \& \, i \neq k \}.
\nonumber
\eea
Here we would like to note that the sign of $u_k$ could be either
positive or negative.
Then we have
\beastar
N_j(u) & = & \sum_{l=0}^{n-j} \left(\int_{R^n_{au,j+l;j}} g(x) dx
- \int_{R^n_{au,j+l+1;j}} g(x) dx\right) \\
&{}& - \sum_{l=0}^{n-j} \Big(\int_{Q^n_{a;j,j+l}\setminus Q^n_{a;j}}
\zeta_{j+l}(x^{j+l-1}) g(x^{n-1}_{a;j+\ell})dx \\
&{}& \qquad \qquad -
\int_{Q^n_{a;j}\setminus Q^n_{a;j,j+l}} \zeta_{j+l}(x^{j+l-1})
g(x^{n-1}_{a;j+\ell})dx\Big) \\
& = & \sum_{l=0}^{n-j} \left(\int_{R^n_{au,j+l;j}\setminus
R^n_{au,j+l+1;j}} g(x) dx
- \int_{R^n_{au,j+l+1;j}\setminus R^n_{au,j+l;j}} g(x) dx
\right) \\
&{}& - \sum_{l=0}^{n-j} \Big(\int_{Q^n_{a;j,j+l}\setminus Q^n_{a;j}}
\zeta_{j+l}(x^{j+l-1}) g(x^{n-1}_{a;j+\ell})dx \\
&{}& \qquad \qquad -
\int_{Q^n_{a;j}\setminus Q^n_{a;j,j+l}}  \zeta_{j+l}(x^{j+l-1})
g(x^{n-1}_{a;j+\ell})dx\Big).
\eeastar

For the simplicity of notations, we will just denote
$$
\zeta_{j+l} =  \zeta_{j+l}(x^{j+l-1})
$$
for the rest of the paper. We now estimate
\beastar
&{}& \left(\int_{R^n_{au,j+l;j}\setminus
R^n_{au,j+l+1;j}} g(x) dx
- \int_{R^n_{au,j+l+1;j}\setminus R^n_{au,j+l;j}} g(x) dx
\right)\\
&{}& \quad -  \left(\int_{Q^n_{a;j,j+l}\setminus Q^n_{a;j} }
 \zeta_{j+l}\cdot g(x^{n-1}_{a;j+\ell})dx -
\int_{Q^n_{a;j}\setminus Q^n_{a;j,j+l}}  \zeta_{j+l}\cdot
g(x^{n-1}_{a;j+\ell})dx\right)
\eeastar
for each $l =0, 1, \cdots, n-j$. We further rewrite it as
\bea
&{}& \left(\int_{R^n_{au,j+l;j}\setminus
R^n_{au,j+l+1;j}} g(x) dx
- \int_{R^n_{au,j+l+1;j}\setminus R^n_{au,j+l;j}} g(x) dx
\right.\nonumber \\
&{}& \quad \left.- \int_{R^n_{au,j+l;j}\setminus
R^n_{au,j+l+1;j}} \zeta_{j+l}\cdot g(x^{n-1}_{a;j+l}) dx
\right.\nonumber\\
&{}& \hskip1.0in \left.- \int_{R^n_{au,j+l+1;j}\setminus R^n_{au,j+l;j}}
 \zeta_{j+l}\cdot g(x^{n-1}_{a;j+l}) dx
\right) \label{eq:first} \\
&{}& \quad + \left(\int_{R^n_{au,j+l;j}\setminus
R^n_{au,j+l+1;j}} \zeta_{j+l}\cdot g(x^{n-1}_{a;j+l}) dx
\right.\nonumber\\
&{}& \hskip1.0in \left.- \int_{R^n_{au,j+l+1;j}
\setminus R^n_{au,j+l;j}}
 \zeta_{j+l}\cdot g(x^{n-1}_{a;j+l}) dx
\right. \nonumber \\
&{}& \quad \left.- \int_{Q^n_{a;j,j+l}\setminus Q^n_{a;j} }
 \zeta_{j+l}\cdot g(x^{n-1}_{a;j+l})dx -
\int_{Q^n_{a;j}\setminus Q^n_{a;j,j+l}}
 \zeta_{j+l}\cdot g(x^{n-1}_{a;j+l})dx\right).
\label{eq:second}
\eea

We estimate (\ref{eq:first}) and (\ref{eq:second}) separately.
We start with (\ref{eq:first}).

The terms in (\ref{eq:first}) can be combined into
\bea\label{eq:newfirst}
& {}& \left(\int_{R^n_{au,j+l;j}\setminus
R^n_{au,j+l+1;j}} g(x) dx
- \int_{R^n_{au,j+l+1;j}\setminus R^n_{au,j+l;j}} g(x) dx
 \right) \nonumber \\
&{}& \quad - \left(\int_{R^n_{au,j+l;j}\setminus
R^n_{au,j+l+1;j}}  \zeta_{j+l}\cdot g(x^{n-1}_{a;j+l}) dx
\right.\nonumber \\
&{}& \hskip1.0in \left.- \int_{R^n_{au,j+l;j}\setminus R^n_{au,j+l;j}}
 \zeta_{j+l}\cdot g(x^{n-1}_{a;j+l}) dx
\right) \nonumber \\
& = & \int_{R^n_{au,j+l;j}\setminus
R^n_{au,j+l+1;j}}\left(g(x)- \zeta_{j+l}\cdot g(x^{n-1}_{a;j+l})\right) dx
\nonumber \\
&{}& \hskip1.0in - \int_{R^n_{au,j+l+1;j}\setminus R^n_{au,j+l;j}} \left(g(x)
- \zeta_{j+l}\cdot g(x^{n-1}_{a;j+l})\right) dx
\eea
Then we obtain the inequality
\beastar
m(R^n_{au,j+l;j}\setminus R^n_{au,j+l+1;j}) \leq
\left(\prod_{i=j}^{j+l} a_i\right)\cdot |\zeta_{j+l+1}u_{j+l+1}|\cdot
\left(\prod_{i=j+2+l}^n |a_i+\zeta_i u_i(a)|\right) \leq |\zeta\cdot u|
\eeastar
and the same for $m(R^n_{au,j+l+1;j}\setminus R^n_{au,j+l;j})$.

We have the bound for the first term of (\ref{eq:newfirst})
\bea
&{}&\Big|\int_{R^n_{au,j+l;j}\setminus
R^n_{au,j+l+1;j}}\left(g(x) -  \zeta_{j+l}\cdot
g(x^{n-1}_{a;j+l})\right) dx\Big|\nonumber\\
& \leq &
\Big|\int_{R^n_{au,j+l;j}\setminus
R^n_{au,j+l+1;j}}\left(g(x) -
g(x^{n-1}_{a;j+l})\right) dx\Big|\label{eq:first1}\\
&{}& +
\int_{R^n_{au,j+l;j}\setminus
R^n_{au,j+l+1;j}}|1 -  \zeta_{j+l}|
g(x^{n-1}_{a;j+l}) dx. \label{eq:first2}
\eea

To get a bound for (\ref{eq:first1}),
we now define
\be\label{eq:Rn-1ajk}
R^{n-1}_{au;j,k}(a) = \{(x_1, \cdots, \widehat x_k, x_{k+1},\cdots, x_n)
\mid (x_1, \cdots, a_k, x_{k+1},\cdots, x_n) \in R^n_{au;j} \}
\ee
for each $j \leq k \leq n$.
Fubini's theorem then implies that the term (\ref{eq:first1})
can be bounded by
\beastar
\int_{R^n_{au,j+l;j}\setminus
R^n_{au,j+l+1;j}}|g(x) - g(x^{n-1}_{a;j+l})| dx &{}& \nonumber\\
= \int_{R^{n-1}_{a;j,j+l}(a)}
\int_{M_{j+l+1}^-(a;u)}^{M_{j+l+1}^+(a;u)}
|g(x) - g(x^{n-1}_{a;j+l})|  dx
\eeastar
where we define the constants
$$
M_k^-(a;\zeta\cdot u) = \min\{ a_k, a_k + \zeta_k u_k(a) \}, \quad
M_k^+(a;\zeta\cdot u) = \max\{ a_k, a_k + \zeta_k u_k(a) \}.
$$
Noting that
$$
M_{j+l+1}^+(a;u) - M_{j+l+1}^-(a;u) = |\zeta_{j+l+1} u_{j+l+1}(a)|,
$$
we have
$$
|x - x^{n-1}_{a;j+l}|\leq |\zeta_{j+l+1} u_{j+l+1}(a)| \leq |\zeta\cdot u|
$$
for all
\be\label{eq:x}
x \in R^n_{au,j+l;j}\setminus
R^n_{au,j+l+1;j} \quad \mbox{ and } \,
\pi_{j+l+1}(x) = \pi_{j+l+1}(x^{n-1}_{a;j+l+1})
\ee
where $\pi_{j+l+1}:\R^n \to \R^{n-1}$ is the coordinate projection
along the $(j+l+1)$-th axis.

Using this preparation, we have the bound for(\ref{eq:first1}) given by
$$
|\zeta\cdot u| \cdot \max |g(x) - g(x^{n-1}_{a;j+l})| :
$$
Here the maximum is taken over all $x, \, x^{n-1}_{a;j+l+1}$ satisfying
(\ref{eq:x}).
\emph{Using the continuity of $g$ and compactness of $Q$}, we have
\be\label{eq:maxM-M+}
|g(x) - g(y)| \leq L_g|x-y|
\ee
for some $L_g > 0$ depending only on $g$. We also note
$$
|a_i + \zeta_i u_i(a)| \leq 1.
$$
Therefore (\ref{eq:first1}) is bounded by
\be\label{eq:first1<}
|\zeta_{j+l} u_{j+l}|L_g |\zeta\cdot u| \leq 2 L_g |\zeta\cdot u||\zeta\cdot u|.
\ee

On the other hand, Fubini's theorem implies that the term (\ref{eq:first2})
can be bounded by
\bea\label{eq:1-zeta}
\int_{R^n_{au,j+l;j}\setminus
R^n_{au,j+l+1;j}}|1 -  \zeta_{j+l}|
g(x^{n-1}_{a;j+l}) dx &{}& \nonumber\\
= \int_{R^{n-1}_{a;j,j+l}(a)}
\int_{M_{j+l+1}^-(a;u)}^{M_{j+l+1}^+(a;u)} |1- \zeta_{j+l}|
g(x^{n-1}_{a;j+l}) dx.
\eea

Now we define $\e_1(\zeta)$ to be
$$
\e_1(\zeta) = \max_{1\leq k \leq n}\left(
\int_{Q^{k-1}}|1 -  \zeta_k(x^{k-1})| dx^{k-1} \right)
$$
as in (\ref{eq:e1zeta}). With this definition and by Fubini's
theorem, we estimate
(\ref{eq:1-zeta}) by
\bea\label{eq:<e1g}
\quad &{}& \int_{R^{n-1}_{a;j,j+l}(a)}
\int_{M_{j+l+1}^-(a;u)}^{M_{j+l+1}^+(a;u)} |1- \zeta_{j+l+1}|
g(x^{n-1}_{a;j+l+1}) dx  \nonumber\\
&\leq & \max g \int_{Q^{j+l}}\int_{Q^{n-(j+l+1)}_{a;j+l}}
\int_{M_{j+l+1}^-(a;u)}^{M_{j+l+1}^+(a;u)} |1- \zeta_{j+l+1}(x^{j+l})|
dx \nonumber \\
& \leq &  \max g |\zeta\cdot u| \int_{Q^{j+l}_{a;j+l}}|1- \zeta_{j+l+1}|
dx \nonumber \\
& \leq & \e_1(\zeta)|\zeta\cdot u| \max g.
\eea
Here we denote
$$
Q^{n-k}_{a;k} = \{(a_k, \widetilde x_k) \mid \widetilde x_k \in Q^{n-(k+1)} \}.
$$
Since $g$ is fixed, we can make $\e_1(\zeta)$ as small as
we want by choosing the cut-off function $\zeta$ as close to 1
as possible in the $L^1$ sense.

Applying the above discussion of (\ref{eq:first1}) and (\ref{eq:first2}), we derive
$$
\Big|\int_{R^n_{au,j+l;j}\setminus
R^n_{au,j+l+1;j}}\left(g(x) - \zeta_{j+l+1}\cdot g(x^{n-1}_{a;j+l+1})\right) dx\Big|
\leq (\e_1(\zeta)\max g + 2 L_g|\zeta\cdot u|)|\zeta\cdot u|
$$
from (\ref{eq:first}) and (\ref{eq:<e1g}).
Here we use the assumption that $|\zeta\cdot u| \leq 2$ and the fact $a_i < 1$.
Similar estimate gives the bound
\be
\Big|\int_{R^n_{au,j+l+1;j} \setminus R^n_{au,j+l;j}}
\left(g(x)-\zeta_{j+l+1}\cdot g(x^{n-1}_{a;j+l+1})\right) dx\Big|
\leq (\e_1(\zeta)\max g + 2L_g|\zeta\cdot u|) |\zeta \cdot u|
\ee
and hence (\ref{eq:first}) is bounded by
\be\label{eq:e22}
2(\e_1(\zeta)\max g+ 2 L_g|\zeta\cdot u|) |\zeta\cdot u|.
\ee

Next, we turn to (\ref{eq:second}).
The four terms in (\ref{eq:second}) can be combined into
\beastar
&{}& \left(\int_{R^n_{au,j+l;j}\setminus
R^n_{au,j+l+1;j}} \zeta_{j+l}\cdot g(x^{n-1}_{a;j+l}) dx
- \int_{R^n_{au,j+l+1;j}\setminus R^n_{au,j+l;j}}
\zeta_{j+l}\cdot g(x^{n-1}_{a;j+l}) dx \right.\nonumber\\
&{}& \quad \left. - \int_{Q^n_{a;j,j+l}\setminus Q^n_{a;j} }
 \zeta_{j+l}\cdot g(x^{n-1}_{a;j+l})dx-
\int_{Q^n_{a;j}\setminus Q^n_{a;j,j+l}}  \zeta_{j+l}\cdot g(x^{n-1}_{a;j+l})dx\right)
\nonumber \\
&{}& = \left(\int_{R^n_{au,j+l;j}\setminus
R^n_{au,j+l+1;j}} \zeta_{j+l}\cdot g(x^{n-1}_{a;j+l}) dx
- \int_{Q^n_{a;j,j+l}\setminus Q^n_{a;j} }  \zeta_{j+l}\cdot g(x^{n-1}_{a;j+l})dx\right)
\nonumber\\
&{}& \quad  -\left(\int_{R^n_{au,j+l+1;j}\setminus R^n_{au,j+l;j}}
 \zeta_{j+l}\cdot g(x^{n-1}_{a;j+l}) dx - \int_{Q^n_{a;j}\setminus Q^n_{a;j,j+l}}
 \zeta_{j+l}\cdot g(x^{n-1}_{a;j+l})dx\right)
\eeastar
Here we can write
\beastar
&{}& \int_{R^n_{au,j+l;j}\setminus
R^n_{au,j+l+1;j}} \zeta_{j+l}\cdot g(x^{n-1}_{a;j+l}) dx
- \int_{Q^n_{a;j,j+l}\setminus Q^n_{a;j} } \zeta_{j+l}\cdot g(x^{n-1}_{a;j+l})dx\\
& = &\int_{\gamma_{j+l+1}(\varphi_{j+l}(Q^n_{a;j})\setminus Q^n_{a;j}))}
 \zeta_{j+l}\cdot g(x^{n-1}_{a;j+l}) dx
- \int_{Q^n_{a;j,j+l}\setminus Q^n_{a;j}} \zeta_{j+l}\cdot g(x^{n-1}_{a;j+l})dx\\
& = & \int_{(\gamma_{j+l+1}(\varphi_{j+l}(Q^n_{a;j})\setminus Q^n_{a;j}))
\setminus (Q^n_{a;j,j+l}\setminus Q^n_{a;j})}
 \zeta_{j+l}\cdot g(x^{n-1}_{a;j+l}) dx\\
&{}& \quad - \int_{(Q^n_{a;j,j+l}\setminus Q^n_{a;j})\setminus
(\gamma_{j+l+1}(\varphi_{j+l}(Q^n_{a;j})\setminus Q^n_{a;j}))}
\zeta_{j+l}\cdot g(x^{n-1}_{a;j+l})dx
\eeastar
where $\gamma_{j+l+1}:Q^n \to Q^n$ is the diffeomorphism
$$
\gamma_{j+l+1} = \varphi_n \circ \cdots \circ \varphi_{j+l+1}.
$$

Now we prove the following lemma

\begin{lem} Let $x \in (\gamma_{j+l+1}(\varphi_{j+l}(Q^n_{a;j,j+l})\setminus Q^n_{a;j}))
\setminus (Q^n_{a;j,j+l}\setminus Q^n_{a;j})$.
For each $0 \leq l < n-j$, $x$ must satisfy
\beastar
0 & \leq & x_i \leq 1 \hskip0.85in \mbox{for } \, 1 \leq i \leq j-1 \\
0 & \leq & x_i \leq a_i \hskip0.8in \mbox{for } \, j \leq i \leq j+l-1 \\
M_k^-(a;\zeta \cdot u) &\leq& x_k \leq M_k^+(a;\zeta\cdot u) \quad \mbox{for some }
j+l+1 < k \leq n\\
a_{j+l} + u_{j+l}(a_{j+l},\widetilde x_{j+l}) & \leq & x_{j+l} \leq
a_{j+l} + \zeta_{j+l}u_{j+l}(a_{j+l},\widetilde x_{j+l}).
\eeastar
For $l = n-j$, we have
$$
(\gamma_{j+l+1}(\varphi_{j+l}(Q^n_{a;j,j+l})\setminus Q^n_{a;j}))
\setminus (Q^n_{a;j,j+l}\setminus Q^n_{a;j}) = \varphi_n(Q^n_{a;j,j+l})\setminus Q^n_{a;j}))
\setminus (Q^n_{a;j,j+l}\setminus Q^n_{a;j})
$$
and
\beastar
0 & \leq & x_i \leq 1 \hskip 0.7in \mbox{for } \, 1 \leq i \leq n-1 \\
a_n + u_n(a_n) & \leq & x_n \leq
a_n + \zeta_n(x^{j+l-1}) u_n(a_n).
\eeastar
\end{lem}
\begin{proof} First consider the case $0 \leq l < n-j$. Note that
\beastar
\gamma_{j+l+1}(x_1, \cdots,x_n) & = &
(x_1, \cdots, x_{j+l}, x_{j+l+1} +
\zeta_{j+l+1}(x^{j+l})u_{j+l+1}(x_{j+l+1},\widetilde x^{j+l+1}),\\
&{}& \quad \cdots x_n + \zeta_n(\varphi_{n-1}\circ\cdots\circ \varphi_{j+l}(x))u_n(x_n)).
\eeastar
Let
\beastar
x & \in & (\gamma_{j+l+1}(\varphi_{j+l}(Q^n_{a;j,j+l})\setminus Q^n_{a;j}))
\setminus (Q^n_{a;j,j+l}\setminus Q^n_{a;j}) \\
&{}& \qquad =
\gamma_{j+l+1}(\varphi_{j+l}(Q^n_{a;j,j+l}))
\setminus (Q^n_{a;j,j+l}\cup Q^n_{a;j}).
\eeastar
Then we first have
$$
a_{j+l} \leq x_{j+l} \leq
a_{j+l} + \zeta_{j+l}u_{j+l}(a_{j+l}, \widetilde x_{j+l}).
$$
(Here if $u_{j+l}(a_{j+l}, \widetilde x_{j+l})< 0$, this inequality is vacuous, i.e.,
no such $x$ exists.)

And for some $k \geq j+l +1$, we have $x_k \geq a_k$ and can write
$$
x_k = y_k + \zeta_k u_k(y_k, \widetilde y_k)
$$
for some $y$ satisfying
$$
y_k \leq a_k.
$$
In particular, we obtain
$$
a_k \leq x_k \leq a_k + \zeta_k u_k(a_k, \widetilde x_k)
$$
and hence we have obtained
$$
(M_k^-(a;\zeta \cdot u) \leq) a_k \leq x_k \leq M_k^+(a;\zeta\cdot u).
$$
The proof of other inequalities are easy and so omitted.

For the case $l = n-j$, we just note
$$
\gamma_{n+1} =id
$$
and then the rest follows.
\end{proof}

Noting $0 \leq a_i \leq 1$, $|\zeta_{j+l}| \leq 1+\e \leq 2$ and
$$
M_k^+(a;\zeta\cdot u) - M_k^-(a;\zeta\cdot u) =|\zeta_k u_k(a)|
$$
we obtain
\bea
&{}& m\left((R^n_{au,j+l;j} \setminus
R^n_{au,j+l+1;j}) \setminus (Q^n_{a;j,j+l}\setminus Q^n_{a;j})\right)
\nonumber\\
& \leq  &
\left(\prod_{j \leq i < j+l} |a_i|\right)\cdot |\zeta_{j+l+1}u_{j+l+1}|
|(1+\zeta_{j+l+1})u_{j+l+1}| \nonumber \\
& \leq & (2+\e)|\zeta\cdot u||u|
\label{eq:1leqln-j}
\eea
for $j \leq l < n-j$. Therefore for  $j \leq l < n-j$, we have proved
\be\label{eq:mRQ}
m\left((R^n_{au,j+l;j}\setminus
R^n_{au,j+l+1;j}) \setminus (Q^n_{a;j,j+l}\setminus Q^n_{a;j})
\right)
\leq 3|\zeta\cdot u||u|.
\ee
Similarly we prove
\be\label{eq:mQR}
m\left((Q^n_{a;j,j+l}\setminus Q^n_{a;j}) \setminus (R^n_{au,j+l;j}\setminus
R^n_{au,j+l+1;j})\right) \leq 3|\zeta\cdot u||u|.
\ee
Then we have proved that the absolute value of (\ref{eq:second})
is less than or equal to
\be\label{eq:abssecond}
6(\max g) |\zeta\cdot u||u|
\ee
for $j \leq l < n-j$. (Here we need to treat the case of $j=1$ slightly
differently but again the same inequality can be shown to hold whose details we
leave for the readers.)

On the other hand, when $l = n-j$, (\ref{eq:second}) can be estimated as
\bea\label{eq:secondn}
&{}&\quad \Big|\int_{(\varphi_n(Q^n_{a;j})\setminus Q^n_{a;j})
\setminus (Q^n_{a;j,n}\setminus Q^n_{a;j})}
\zeta_{j+l}\cdot g(x^{n-1}_{a;j+l}) dx \nonumber\\
&{}&- \int_{(Q^n_{a;j,n}\setminus Q^n_{a;j})\setminus
(\varphi_n(Q^n_{a;j})\setminus Q^n_{a;j})}
\zeta_{j+l}\cdot g(x^{n-1}_{a;j+l})dx\Big| \nonumber \\
& \leq & 2 \int_{Q^{n-1}} \int_{M^-_n(a;\zeta\cdot u)}^{M^+_n(a;\zeta\cdot u)}
\zeta_n(x^{n-1})g(x^{n-1}_{a;n})\, dx_n dx^{n-1}\nonumber \\
& \leq & 4 |\zeta_n u_n(a_n)| \int_{Q^{n-1}}
|1 - \zeta(x^{n-1})| g(x^{n-1}_{a;n})\, dx^{n-1}\nonumber \\
& \leq & 4 \e_1(\zeta)\max g|\zeta\cdot u|.
\eea
Therefore combining (\ref{eq:abssecond}) and (\ref{eq:secondn}),
we have proved that (\ref{eq:second}) is bounded by
\be\label{eq:secondall}
6 (\max g) |\zeta \cdot u||u| + 4 \e_1(\zeta)\max g|\zeta \cdot u|
\ee
for any $1 \leq j \leq n$.

Combining (\ref{eq:secondall}) and (\ref{eq:e22}), we have finally
obtained
$$
|N(u)| \leq \left(8 \max g|\zeta\cdot u| + 6 \max g \e_1(\zeta)
+ 4 L_4|\zeta\cdot u| \right)|u|.
$$
This finishes the proof of Lemma \ref{lemma-N(u)} by setting
$C_4 = \max\{8\max g, 4\}$.


\begin{thebibliography}{Mu2}

\bibitem[Br]{brenier:jams} Brenier, Y., {\it The least action principle and the
related concept of generalized flows for incompressible perfect fluids},
J. Amer. Math. Soc. {\bf 2} (1989), 225-255.


\bibitem[DM]{DM} Dacorogna, B. and Moser, J., {\it On a partial differential
equation involving the Jacobian determinant}, Ann. Inst. Henri Poincar\'e,
{\bf 7} (1990), 1-26.

\bibitem[D]{donald} Donaldson, S. K.,
{\it An application of gauge theory to four-dimensional topology},
J. Differential Geom. {\bf 18} (1983), no. 2, 279--315.

\bibitem[E]{eliash:rigid} Eliashberg, Y., {\it A theorem on the structure
of wave fronts and applications in symplectic topology},
Funct. Anal. and its Appl. {\bf 21} (1987), 227-232.

\bibitem[F]{fathi} Fathi, A., {\it Structure of the group of
homeomorphisms preserving a good measure on a compact manifold},
Ann. Scient. \`Ec. Norm. Sup. {\bf 13} (1980), 45-93.

\bibitem[G1]{gromov:pseudo} Gromov, M., {\it Pseudo-holomorphic curves in
symplectic manifolds}, Invent. Math. {\bf 82} (1985), 307-347.

\bibitem[G2]{gromov:metric} Gromov, M., Metric Structures for Riemannian
and Non-Riemannian Spaces, Progress in Math. 152, Birkh\"auser, Boston,
1999.

\bibitem[La]{lang} Lang, S. Real Analysis, 2nd edition, Addison Wesley,
1983, Reading, MA.

\bibitem[Mo]{moser} Moser, J., {\it On the volume elements on a manifold},
Trans. Amer. Math. Soc. {\bf 120} (1965), 286-294.

\bibitem[Mu1]{munkres:S2} Munkres, J., {\it Differentiable isotopies on the
2-sphere}, Michigan Math. J. {\bf 7} (1960), 193-197.

\bibitem[Mu2]{munkres:smoothing} Munkres, J., {\it Obstructions to the smoothing
of piecewise-differentiable homeomorphisms}, Ann. Math. {\bf 72} (1960), 521-554.


\bibitem[OF]{oh:icm} Oh, Y.-G., Fukaya, K., {\it
Floer homology in symplectic geometry and in mirror symmetry},
Proceedings of ICM-2006, 879 - 905, Madrid, Spain, 2006.

\bibitem[OM]{oh:hameo1} Oh, Y.-G., M\"uller, S., {\it The group of Hamiltonian homeomorphisms
and $C^0$ symplectic topology}, J. Symp. Geom. (to appear), arXiv:math.SG/0402210.

\bibitem[OU]{oxtoby} Oxtoby, J. C. and Ulam, S. M., {\it Measure preserving
homeomorphisms and metrical transitivity},
Ann. Math. {\bf 42} (1941), 874-920.

\bibitem[Sh]{shnirel:gafa} Shnirelman, A., {\it Generalized fluid flows,
their approximation and applications,} Geom. Funct. Anal. {\bf 4} (1994),
586 - 620.

\bibitem[Vi]{villani} Villani, C., Topics in Optimal Transportation, Graduate Studies in Math.
{\bf 58}, AMS, 2003, Providence.
\end{thebibliography}
\end{document}